\documentclass[10pt]{amsart}
\usepackage{hyperref}
\usepackage{amssymb,amsmath}
\newtheorem{theorem}{Theorem}[section]
\newtheorem{proposition}{Proposition}[section]
\newtheorem{lemma}{Lemma}[section]
\newtheorem{corollary}{Corollary}[section]
\newtheorem{remark}{Remark}[section]

\newtheorem{definition}{Definition}[section]
\numberwithin{equation}{section}
\newcommand{\bee}{\begin{equation}}
\newcommand{\eee}{\end{equation}}
\begin{document}
\title[]
{On the Dynamics of a Degenerate Parabolic Equation: Global Bifurcation of Stationary States and Convergence}
\author[]{Nikos I. Karachalios $^{\dag}$
and Nikos B. Zographopoulos$^{\ddag}$}
\thanks
{$^{\dag}$ Department of Mathematics,
University of the Aegean, GR 83200 Karlovassi, Samos, GREECE (karan@aegean.gr)}
\thanks
{$^{\ddag}$ Department of Applied Mathematics,
University of Crete,
71409 Heraklion, Crete, GREECE (nzogr@tem.uoc.gr).}
%
%
%
\centerline{\today}
\begin{abstract}
We study the dynamics of a degenerate parabolic equation with a
variable, generally non-smooth diffusion coefficient, which may
vanish at some points  or be unbounded.  We show the
existence of a global branch of nonnegative stationary states,
covering both the cases of a bounded and an unbounded domain.
The global bifurcation of stationary states, implies-in
conjuction with the definition of a gradient dynamical system in
the natural phase space-that at
least in the case of a bounded domain, any solution with
nonnegative initial data tends to the trivial or the nonnegative
equilibrium.  Applications of the global bifurcation result to
general degenerate semilinear as well as to quasilinear elliptic
equations, are also discussed.
\end{abstract}
\maketitle
%
%
%
%
\section{Introduction}
The mathematical modelling of various physical processes, where
spatial heterogeneity has a primary role, has usually as a result,
the derivation of nonlinear evolution equations with variable
diffusion, or dispersion. Applications are ranging from physics to
biology. To name but a few, equations of such a type have been
successfully applied to the heat propagation in heterogeneous
materials \cite{DL,JRocha,kamin2,kamin3}, the study of transport
of electron temperature in a confined plasma \cite{kamina}, the
propagation of varying amplitude waves in a nonlinear medium
\cite{sulem} (and \cite{Craig} for linear  Schr\"{o}dinger
equation), to the study of electromagnetic phenomena in
nonhomogeneous superconductors \cite{ChapRich97,
JKR92,JianSong01,JimMor99} and the dynamics of Josephson junctions
\cite{Flytzanis, Flytzanis1}, to epidemiology and the growth and
control of brain tumors \cite{mur03}.

In this work we continue the study, initiated in \cite{kz??}, of
the qualitative behavior of solutions of some degenerate evolution equations (involving degenerate coefficients). Work \cite{kz??} concerns the asymptotic behavior of solutions, of a complex evolution equation of Ginzburg-Landau type. Here we study 
the following semilinear parabolic equation with variable, nonnegative
diffusion coefficient, defined on an arbitrary domain (bounded or
unbounded)\ $\Omega \subseteq \mathbb{R}^N$, $N\geq 2$,\
\begin{eqnarray} \label{eq1.0}
\partial_t \phi - \mathrm{div}(\sigma(x)\, \nabla \phi ) &-& \lambda\, \phi
+|\phi|^{2\gamma}\phi=0,\;\;x\in\Omega,\;\;t>0, \nonumber\\
\phi(x,0)&=& \phi_0(x),\;\;x\in\Omega,\\
\phi|_{\partial\Omega}&=&0,\;t>0\nonumber.
\end{eqnarray}
Equation (\ref{eq1.0}) can be derived as a simple model for
neutron diffusion (feedback control of nuclear reactor)
\cite{DL,kasten69}. In this case $\phi$ (which must be
nonnegative) and  $\sigma$ stand for the neutron flux and neutron
diffusion respectively.

The degeneracy of problem  (\ref{eq1.0}) is considered in the
sense that the measurable, nonnegative diffusion coefficient
$\sigma$, is allowed to have at most a finite number of
(essential) zeroes, at some points or even to be unbounded. The
point of departure for the consideration of suitable assumptions
on the diffusion coefficient is the work \cite{cm00}, where the
degenerate elliptic problem  is studied: we
assume that the function $\sigma:\Omega\rightarrow\mathbb{R}$
satisfies the following assumptions\vspace{.2cm}
\newline
$(\mathcal{H_{\alpha}})$\ \ $\sigma\in L^1_{\mathrm{loc}}(\Omega )$ and for
some $\alpha\in (0,2)$, $\liminf_{x\rightarrow
z}|x-z|^{-\alpha}\sigma(x)>0$, for every $z
\in\overline{\Omega}$, when the domain $\mathrm{\Omega}$ is bounded,\vspace{.2cm}
\newline
$(\mathcal{H^{\infty}_{\beta}})$\ \ $\sigma$ satisfies condition
$(\mathcal{H_{\alpha}})$ and $\liminf_{|x| \rightarrow
\infty}|x|^{-\beta} \sigma(x)>0$,\ for some $\beta>2$, when the
domain $\mathrm{\Omega}$ is unbounded. \ \vspace{.2cm}

The assumptions\ $(\mathcal{H_{\alpha}})$\ and\
$(\mathcal{H^{\infty}_{\beta}})$\ imply (see \cite[Lemma
2.2]{cm00}) that (i) the set of zeroes is finite, (ii) the
function\ $\sigma$\ could be non smooth (cannot be of class\
$C^2$, if $\alpha \in (0,2)$\ and it cannot have bounded
derivatives if\ $\alpha \in (0,1)$).\ Moreover, in the unbounded
domain case the function\ $\sigma$\ has to be unbounded. The
approach in \cite{cm00}, was based on Caffarelli-Kohn-Nirenberg
type inequalities (see (\ref{CKN})). For some recent results
concerning these inequalities and their applications to the study of elliptic
equations, we refer to \cite{abd03, cat01, fel03}. \vspace{0.2cm}

The physical motivation of the assumption
$(\mathcal{H_{\alpha}})$, is related to the modelling of reaction
diffusion processes in composite materials, occupying a bounded
domain $\mathrm{\Omega}$, which at some points they behave as {\em
perfect insulators}. Following \cite[pg. 79]{DL}, when at some
points the medium is perfectly insulating, it is natural to assume
that $\sigma(x)$ vanishes at these points. On the other hand, when
condition $(\mathcal{H^{\infty}_{\beta}})$ is satisfied,  it
follows from  \cite[Lemma 2.2]{cm00}, that in addition, the diffusion coefficient
has to be unbounded.
Physically, this situation corresponds to a nonhomogeneous medium,
occupying the unbounded domain $\mathrm{\Omega}$, which behaves as
a perfect conductor in $\Omega\setminus \mathrm{B}_R(0)$ (see
\cite[pg.79]{DL}), and as a perfect insulator in a finite number
of points in $\mathrm{B}_R(0)$. Note that when
$\partial\Omega\neq\emptyset$, the function $\sigma(x)$, need not
be locally bounded. These conditions arise in various simple
transport models of electron temperature in a confined plasma. See
\cite{kamin3} for a discussion concerning the one-dimensional
case: the electron thermal diffusion is density dependent such
that it vanishes with density, rendering the problem singular.
Note that in various  diffusion processes, the equations
involve  diffusion $\sigma(x)\sim |x|^{\alpha}$, $\alpha < N$:
We refer to \cite{kamina,kaminb} for equations describing heat
propagation.

The main purpose of this work is to combine basic results from the theory of infinite dimensional dynamical systems
and bifurcation theory, to give a description of the dynamics of (\ref{eq1.0}). {\em We remark here the crucial role of the conditions $(\mathcal{H_{\alpha}})$ and $(\mathcal{H^{\infty}_{\beta}})$ on the ``degeneracy exponents'' $\alpha, \beta$ which give rise to necessary compactness properties of  various linear and nonlinear operators} associated to the study of (\ref{eq1.0}) and its related stationary problem (a degenerate elliptic equation).
We are restricted in the case $N\geq 2$ since the case $N=1$, despite its similarities with the higher dimensional case with respect to the definition and properties of the appropriate functional setting, recovers also important differences. For the definition and properties of the related function spaces and detailed discussions on one dimensional versions of generalized Hardy and Caffarelli-Kohn-Nirenberg inequalities, we refer to \cite{cmad,cat01}.

More precisely, the first part of the present work is devoted to some results
concerning the existence of a global attractor. While the result in  \cite{kz??}, for the complex evolution equation, concerns the existence of a global attractor in $L^2(\Omega)$, here it is verified that the dynamical
system associated to (\ref{eq1.0}) is a {\em gradient system}, and
that {\em there exists a connected global attractor} in the
weighted Sobolev space $D^{1,2}_0(\Omega;\sigma)$, the closure of 
$C_{0}^{\infty}(\Omega)$\ with respect to the norm
$||\phi||^2_{D_{0}^{1,2}(\Omega,\sigma)}= \int_{\Omega} \sigma(x)\,
|\nabla \phi|^2$. 
This space appears to be the natural
energy space for (\ref{eq1.0}). The main result of
Section 3, can be stated by the following theorem.
\begin{theorem}
\label{main1} Let\ $\Omega \subseteq \mathbb{R}^N$,\ $N \geq 2$,\
be an arbitrary domain (bounded or unbounded). Assume that\
$\sigma$\ satisfies condition\ $(\mathcal{H_{\alpha}})$ or
$(\mathcal{H^{\infty}_{\beta}})$  and
$$0<\gamma\leq\frac{2-\alpha}{2(N-2+\alpha)}:=\gamma^*.$$ Equation
(\ref{eq1.0}) defines a semiflow
$$\mathcal{S}(t):D^{1,2}_0(\Omega,\sigma)\rightarrow
D^{1,2}_0(\Omega,\sigma),$$ which possesses a global attractor
${\mathcal A}$ in  $D^{1,2}_0(\Omega;\sigma)$. Let $\mathcal{E}$
denote the (bounded) set of equilibrium points of
$\mathcal{S}(t)$. For each positive orbit $\phi$ lying in
$\mathcal{A}$ the limit set $\omega(\phi)$ is a connected subset
of $\mathcal{E}$ on which
$\mathcal{J}:D^{1,2}_0(\Omega;\sigma)\rightarrow\mathbb{R}$, 
\begin{eqnarray*}
\mathcal{J}(\phi):=\frac{1}{2}\int_{\Omega}\sigma(x)|\nabla\phi|^2dx-\frac{\lambda}{2}\int_{\Omega}|\phi|^2dx+\frac{1}{2\gamma+2}\int_{\Omega}|\phi|^{2\gamma +2}dx,
\end{eqnarray*}
the Lyapunov functional associated to $\mathcal{S}(t)$, is constant.
If $\mathcal{E}$ is totally disconnected (in particular if
$\mathcal{E}$ is countable), the limit
$$z_{+}=\lim_{t\rightarrow+\infty}\phi (t),$$
exists and is an equilibrium point. Furthermore, any solution of
(\ref{eq1.0}), tends to an equilibrium point as
$t\rightarrow+\infty$.
\end{theorem}

Further analysis is carried out, regarding the  bifurcation of the
corresponding steady states with respect to the parameter
$\lambda\in\mathbb{R}$. More precisely, we prove the existence of
{\em a global branch of nonnegative solutions} for the  equation
\begin{equation} \label{eq1.1}
\begin{array}{ll}
-\mathrm{div}(\sigma(x) \nabla u) = \lambda\, u - |u|^{2\gamma}\,
u,\;\;\mbox{in}\;\;\Omega, \vspace{0.2cm} \\
\;\;\;\;\;\;\;\;\;\;\;\;\;\;\;\, u |_{\partial \Omega} = 0,
\end{array}
\end{equation}
bifurcating from the trivial solution at\ $(\lambda_1 ,0)$,\
where\ $\lambda_1$, is the positive principal eigenvalue of the
corresponding linear problem
\begin{equation} \label{eq1.2}
\begin{array}{ll}
-\mathrm{div}(\sigma(x) \nabla u) = \lambda\,
u,\;\;\mbox{in}\;\;\Omega, \vspace{0.2cm} \\
\;\;\;\;\;\;\;\;\;\;\;\;\;\;\;\, u |_{\partial \Omega} = 0.
\end{array}
\end{equation}
This is the main result of Section 4, described by the following Theorem.
\begin{theorem}
\label{main2} Let\ $\Omega \subseteq \mathbb{R}^N$,\ $N \geq 2$,\
be an arbitrary domain (bounded or unbounded). Assume that\ $\sigma$\
satisfies condition\ $(\mathcal{H_{\alpha}})$ or
$(\mathcal{H^{\infty}_{\beta}})$,
and $$0< \gamma < \frac{2-\alpha}{2(N-2+\alpha)}.$$ Then, the principal
eigenvalue\ $\lambda_1$\ of (\ref{eq1.2}) is a bifurcating point
of the problem (\ref{eq1.1}) and\ $C_{\lambda_1}$\ is a global
branch of nonnegative solutions, which "bends" to the right of\
$\lambda_1$. For any fixed $\lambda >\lambda_1$ these solutions
are unique.
\end{theorem}

The technique leading to the global bifurcation result, is included
in the general strategy of the approximation of solutions of a
degenerate partial differential equation, by constructing an
approximate sequence  of solutions of nondegenerate problems. The
approximation procedure has been successfully applied to evolution
\cite{EMR91, Feiresl}, and to stationary problems \cite{dan91,
dan98, ors02}, and in the context of bifurcation theory
\cite{amb97, bl04, eg00, giac98}.

One of the main difficulties  arising, on the attempt to establish
the global character of the branch of nonnegative solutions for
(\ref{eq1.1}), is {\em that Harnack-type Inequalities are not valid in
general} (see \cite[Remark 3.2]{dkn97}). This is a common fact for
non-uniformly elliptic equations \cite{giltru77}. However, we refer
to  \cite{abd03} and the references therein for generalized
Harnack-type inequalities, applied to degenerate elliptic
equations. Distinguishing between the bounded and the unbounded
domain case, we consider different families of approximate
boundary value problems.

When $\Omega$ is bounded, \cite[Lemma 2.2 and Remark 2.3]{cm00} implies that under assumption $(\mathcal{H_{\alpha}})$,  there exists a finite set $Z:=\{z_1,\ldots ,z_k\}\subset\overline{\Omega}$ and $r,\delta>0$, such that the balls of center $z_i$ and radius $r$, $B_r(z_i)$, $i=1,\ldots ,k$, {\em are pairwise disjoint} and
\begin{equation}
\begin{array}{l}
\label{charw}
(i)\;\;\;\sigma (x)\geq \delta |x-z_i|^{\alpha}\;\;\mbox{for}\;\;x\in B_r(z_i)\cap\Omega,\;\;i=1,\ldots ,k,\\
(ii)\;\;\sigma (x)\geq\delta,\;\mbox{for}\;\;x\in\Omega\setminus\bigcup_i B_r(z_i).
\end{array}
\end{equation}
Moreover if $\sigma$ satisfies $(\mathcal{H_{\alpha}})$, then $\sigma\geq 0$ in $\overline{\Omega}$, the set of zeroes of $\sigma$  $Z_{\sigma}:=\{z\in\overline{\Omega}:\sigma(z)=0\}$ is finite, and $Z_{\sigma}\subseteq Z$ (\cite[Remark 2.3]{cm00}). It is not a loss of a generality to assume that $Z_{\sigma}=Z$.

For convenience and simplicity,  in the bounded domain case, we consider as a model for the diffusion coefficient, the function
\begin{eqnarray}
\label{siga}
\sigma(x)=|x|^{\alpha},\;\;\alpha\in (0,2),
\end{eqnarray}
satisfying $(\mathcal{H_{\alpha}})$.   Quite naturally, we construct a family of approximating {\em nondegenerate} problems as follows:
Setting $\Omega_r:=\Omega\setminus B_r(0)$, we consider the boundary value problems
\[
(P)_r\;\;\;\;\;\; \biggl\{
\begin{array}{ll}
-\mathrm{div}(\sigma (x) \nabla u) = \lambda\, u - |u|^{2\gamma}\, u,\;\;\mbox{in}\;\;\Omega_r,\vspace{0.2cm} \\
\;\;\;\;\;\;\;\;\;\;\;\;\;\;u |_{\partial \Omega_r} = 0.
\end{array}
\]
From the characterization (\ref{charw}), problems $(P)_r$ are
non-degenerate, and it can be shown that for fixed $r>0$, there
exists a global branch of positive solutions (see Definition
\ref{def1.1}), by using Harnack type inequalities. The next step
is to prove that the limit of the approximating family $(P)_r$, as
$r\rightarrow 0$, preserves the same property, thus Theorem
\ref{main2}.

When $\Omega$ is unbounded \cite[Lemma 2.2 and Remark 2.3]{cm00}
implies that under $(\mathcal{H_{\beta}^{\infty}})$, in addition
to (\ref{charw}), there exists $R>0$, such that
$\overline{\mathrm{B}_r(z_i)}\subset\mathrm{B}_R(0)$ for every
$i,\ldots,k$ and
\begin{equation}
\begin{array}{l}
\label{charw1}
(iii)\;\;\sigma(x)\geq \delta |x|^{\beta},\;\;\mbox{for}\;\;x\in\Omega\;\;|x|>R.
\end{array}
\end{equation}
In the unbounded domain case we consider as a model, the diffusion coefficient
\begin{eqnarray}
\label{sigab}
\sigma(x)=|x|^{\alpha}+|x|^{\beta},\;\;\alpha\in (0,2),\;\;\beta>2,
\end{eqnarray}
satisfying $(\mathcal{H_{\beta}^{\infty}})$. Note that since
$\sigma$ is unbounded the Harnack inequality is still not
applicable. To approximate (\ref{eq1.1}) defined in the unbounded
domain ($\Omega\subseteq\mathbb{R}^N$), this time we consider the
approximate family of boundary value problems in
$\Omega_R:=\Omega\cap B_{R}(0)$:\
\[
(P)_R\;\;\;\;\;\; \biggl\{
\begin{array}{ll}
-\mathrm{div}(\sigma(x) \nabla u) = \lambda\, u
-|u|^{2\gamma}\, u,\;\;\mbox{in}\;\;\Omega_R, \vspace{0.2cm} \\
\;\;\;\;\;\;\;\;\;\;\;\;\;\,u |_{\partial
\Omega_R} = 0.
\end{array}
\]
Theorem \ref{main2}  holds for $(P)_R$ and the claim is that as
$R\rightarrow\infty$ the theorem remains valid at the limit.

To establish the properties of the
principal eigenvalues corresponding to  both of the approximating
problems $(P)_r$ and $(P)_R$, we prefer an alternative proof,
based on an appropriate adaptation of Picone's Identity. This
identity has been used in \cite{all01, ah98, ah99}, where the
author established certain properties of the principal eigenvalue of
the p-Laplacian operator, and extends Sturm Theorems to degenerate
elliptic equations.

Furthermore, we note that the presented method is
applicable independently of the shape of $\Omega$.\ In general, the
situation becomes more complicated for non-uniformly elliptic
problems in terms of\ $u$. As an example of the appearance of {\em
local bifurcation}, for such a type of equation, we refer to
\cite{zs02}.

A general treatment of degenerate elliptic equations is provided
by the monograph \cite{dkn97}, focusing on the existence and
properties of solutions (the issue of global bifurcation in the
degenerate case is not addressed). Especially in the unbounded
domain case, the problems are non-degenerate (at least in the sense of
degeneracy, imposed by assumption (\ref{charw1})). In \cite{eg00} a
global bifurcation result is proved for a degenerate semilinear
elliptic equation, with a degenerate diffusion coefficient of
"critical exponent'' (inducing non-compactness). Recent global
bifurcation results for non-degenerate problems are included in
the works \cite{arc01, eg00, giac98, gr00, lop99, rab01, rab03}.
For an overview, we also refer to the latest monographs
\cite{bt03, Hans04}.  \vspace{.2cm}
\\
It is our intention to use Theorem \ref{main2} as a  main tool,
for a more detailed description of the asymptotic behavior of
solutions of (\ref{eq1.0}), at least for the case of a bounded domain. A consequence of Theorem \ref{main2}
is that for any fixed  $\lambda >\lambda_1$, the set $\mathcal{E}$
includes the trivial, the unique
nonnegative solution of (\ref{eq1.0}) and its (unique) nonpositive
reflection. A combination of Theorems  \ref{main1}-\ref{main2}
could be used to design an intuitive picture for the dynamics of
(\ref{eq1.0}): It seems that the system undergoes through
$\lambda_1$ {\em a pitchfork bifurcation of supercritical type,
where exchange of stability holds}, i.e., the trivial solution is
stable when\ $\lambda < \lambda_1$,\ while for\ $\lambda_1 <
\lambda$\ the nonnegative (nonpositive) solution of the global
branch become the stable stationary state. Section 5 is devoted to
some remarks related to the rigorous verification of the
bifurcation picture for (\ref{eq1.1}).  The fact
that solutions of (\ref{eq1.0}) with nonnegative initial data,
remain nonnegative for all times (a "maximum principle"
property), and the stability analysis of the unique nonnegative steady state, in conjuction with  \cite[Theorem 2.7]{Ball2}, implies
the following
\setcounter{corollary}{2}
\begin{corollary}
\label{tsibadyoball} Assume that condition $(\mathcal{H_{\alpha}})$, holds. If $\phi_0\geq 0$ a.e in $\Omega$, any
solution $\phi(t)$ of (\ref{eq1.0}), tends to  either the trivial
or the unique nonnegative equilibrium point, as
$t\rightarrow\infty$.
\end{corollary}
As it is expected, the nonnegative steady state is a global minimizer for the Lyapunov functional (Remark \ref{minLyap}). A comment on the role of the ``degeneracy exponent'' $\alpha$ and a discussion concerning some possible further developments with respect to the case of noncompactness, is given in Remarks \ref{degex}, \ref{CrucComments}. 

We conclude by mentioning the main results, on the convergence of
globally defined and bounded solutions of evolution equations to
rest points, as $t\rightarrow\infty$.

For scalar parabolic equations we refer to \cite{Matanoa, Matanob}
and \cite{Zelenyak} for convergence to a single equilibrium.  In
\cite{Lions} the result is proved for a semilinear heat equation
defined in a higher dimensional domain, assuming a special
structure of the set of rest points (semistable solutions). In \cite {JKR92a},
convergence to a unique rest point, at least for the scalar case,
is proved without the hypothesis that the set of rest points is
totally disconnected. The same result is extended to semilinear
parabolic and wave equations considered in multidimensional
domains in\cite{ Jendoubi2, Jendoubi3,Jendoubi1}, when the
nonlinearity is analytic. For a scalar degenerate parabolic
problem (porous medium equation) a positive answer is given in
\cite{Feiresl}. In the recent work \cite{Busca}, the result of
convergence to a (single) equilibrium is extended to a semilinear
parabolic equation in $\mathbb{R}^N$:  The main difficulty in the
unbounded domain case is that even there exists a unique rest
point $z$ (radial with respect to $0$), the $\omega$-limit set, may contains infinite many
distinct translates of  $z$.  The authors introduce a new method,
by defining moments of energy, which can discriminate against
different translates of a rest point. The work \cite{Busca}
provides  also a brief but complete review of the existing results
and methods. For a more detailed survey we refer to
\cite{Polacik}.

In the case of non-autonomous systems or in the case where
uniqueness of solutions of the evolution equation is not expected,
the question on the convergence of solutions to rest points, and
generally, on the existence of a global attractor, is discussed
through the framework of generalized processes and semiflows in
\cite{Ball1a,BallNS,Ball2}. Applications include  nonautonomous
semilinear wave and parabolic equations, or equations involving
non-Lipschitz nonlinearities, and Navier-Stokes equations.  
%
\section{Preliminaries}
{\bf Function spaces and formulation of the problem.}\ \
We recall some of the basic results on functional
spaces defined in \cite{cm00}.
\setcounter{equation}{0}
Let\ $N \geq 2$,\ $\alpha \in (0,2)$\ and
\[
2^{*}_{\alpha} = : \left\{
\begin{array}{ll}
\frac{4}{\alpha}\in (2,+\infty), \;\;\;\; \mbox{if}  \;\; \alpha \in (0,2),\;\;N=2,\\
\;\; \; \frac{2N}{N-2+\alpha}\in \left(2,\frac{2N}{N-2}\right),
\;\; \mbox{if} \;\; \alpha\in (0,2),\;\;N\geq 3.
\end{array}
\right.
\]

The exponent\ $2^{*}_{\alpha}$,\ has the role of the critical
exponent in the classical Sobolev embeddings. The following
Caffarelli-Kohn-Nirenberg inequality holds, for a constant\ $c$\
depending only on\ $\beta, N$,\
\begin{equation} \label{CKN}
\left(\int_{\mathbb{R}^N}|\phi|^{2^{*}_{\alpha}}\, dx
\right)^{\frac{2}{2^{*}_{\alpha}}} \leq
c\int_{\mathbb{R}^N}|x|^{\beta}|\nabla\phi |^2 dx,\;\mbox{for
every}\;\phi\in C^{\infty}_0(\mathbb{R}^N).
\end{equation}

By using (\ref{CKN}) and conditions\ $(\mathcal{H_{\alpha}})$\
and\ $(\mathcal{H^{\infty}_{\beta}})$, it is proved in
\cite[Proposition 2.5]{cm00}, the following generalized version of
(\ref{CKN}),
\begin{equation} \label{CKNg}
\left(\int_{\Omega}|\phi|^{2^{*}_{\alpha}}\, dx
\right)^{\frac{2}{2^{*}_{\alpha}}}\leq K\; \int_{\Omega}\sigma
(x)|\nabla\phi |^2 dx,\;\mbox{for every}\;\phi\in
C^{\infty}_0(\Omega).
\end{equation}

As a consequence of (\ref{CKN}) and (\ref{CKNg}), we have the following generalized version of
Poincar\'{e} inequality (\cite[Corollary 2.6-Proposition 3.5]{cm00}, see also \cite[Section 5]{kz??}).
\begin{lemma}
\label{gen} Let $\Omega$ be a bounded (unbounded) domain of $\mathbb{R}^N$, $N\geq 2$
and assume that condition $(\mathcal{H_{\alpha}})$ ($(\mathcal{H^{\infty}_{\beta}})$) is satisfied.
Then there exists a constant $c>0$, such that
\begin{equation} \label{poi}
\int_{\Omega}|\phi|^2dx\leq c
\int_{\Omega}\sigma(x)|\nabla\phi|^2dx,\;\;\mbox{for
every}\;\;\phi\in\mathrm{C}^{\infty}_0(\Omega).
\end{equation}
\end{lemma}
We emphasize that inequalities (\ref{CKN}),(\ref{CKNg}) (and  (\ref{poi}) in the case of a bounded domain), hold for some $\alpha\in (0,2]$. However,  the case $a=2$ can be considered as a ``critical case'' with respect {\em to compactness of various embeddings}, even in the bounded domain case.  Moreover, condition\ $(\mathcal{H_{\alpha}})$\
is optimal in the following sense:  For\ $\alpha>2$\ there exist
functions such that (\ref{poi}) is not satisfied \cite{cm00}. Note also that in the case of an unbounded domain, (\ref{poi}) does not hold in general, if  $\beta\leq 2$ in $(\mathcal{H^{\infty}_{\beta}})$. We refer also to the
examples of \cite{abd03}.

The natural energy space for the problems (\ref{eq1.0}) and (\ref{eq1.1}) involves the space
$D_{0}^{1,2}(\Omega,\sigma)$, defined as  the closure of\
$C_{0}^{\infty}(\Omega)$\ with respect to the norm
\[
||\phi||_{D_{0}^{1,2}(\Omega,\sigma)} := \biggl ( \int_{\Omega} \sigma(x)\,
|\nabla \phi|^2  \biggr )^{1/2}.
\]
The space\ $D_{0}^{1,2}(\Omega,\sigma)$\ is a Hilbert space with
respect to the scalar product
\begin{eqnarray*}
(\phi,\psi)_{\sigma}:=\int_{\Omega}\sigma(x)\nabla\phi\nabla{\psi}\,dx,\;\;\mbox{for
every}\;\;\phi\;\psi\in D^{1,2}_0(\Omega,\sigma).
\end{eqnarray*}

The following two lemmas refer to the continuous and compact inclusions of
$D^{1,2}_0(\Omega,\sigma)$ \cite[Propositions 3.3-3.5]{cm00}.
\begin{lemma} \label{lemma2.1}
Assume that\ $\Omega$\ is a bounded domain in\ $\mathbb{R}^N$, $N\geq 2$\ and
\ $\sigma$ satisfies\
$(\mathcal{H_{\alpha}})$.\ Then the following embeddings hold:

i)\ \ $D_{0}^{1,2}(\Omega,\sigma) \hookrightarrow
L^{2^{*}_{\alpha}}(\Omega)$\ continuously, \vspace{0.1cm}

ii)\ $D_{0}^{1,2}(\Omega,\sigma) \hookrightarrow L^p (\Omega)$\
compactly if\ $p \in [1,2^{*}_{\alpha})$.
\end{lemma}
\begin{lemma} \label{lemma2.2}
Assume that\ $\Omega$\ is an unbounded domain in\ $\mathbb{R}^N$, $N\geq 2$,
and \ $\sigma$ satisfies
$(\mathcal{H^{\infty}_{\beta}})$.\
Then the following embeddings hold:

i)\ \ $D_{0}^{1,2}(\Omega,\sigma) \hookrightarrow L^{p}(\Omega)$\
continuously for every\ $p \in [2^{*}_{\beta}, 2^{*}_{\alpha}]$,

ii)\ $D_{0}^{1,2}(\Omega,\sigma) \hookrightarrow L^{p}(\Omega)$\
compactly if\ $p \in (2^{*}_{\beta}, 2^{*}_{\alpha})$.
\end{lemma}
\setcounter{remark}{3}
\begin{remark} 
\label{remacom}
It is crucial to note that as a special
case, the embedding $D^{1,2}_0(\Omega,\sigma)\subset
L^2(\Omega)$ is compact if either conditions $(\mathcal{H_{\alpha}})$
or $(\mathcal{H^{\infty}_{\beta}})$ hold:
Observe that\ $\beta>2$\ implies\ $2^{*}_{\beta} =
\frac{2N}{N-2+\beta} <2$, i.e. $2 \in (2^{*}_{\beta}, 2^{*}_{\alpha})$. In the unbounded domain case, we need $\sigma$ to grow {\em faster than quadratically at infinity}, to ensure compactness.
We also stress the fact, that since $\sigma$ is not in
$L^{\infty}_{\mathrm{loc}}(\Omega)$,  there is not in
general any inclusion relation between the space
$D^{1,2}_0(\Omega,\sigma)$ and the standard Sobolev space
$\mathrm{H}^1_0(\Omega)$.
\end{remark}
To justify the natural energy space for (\ref{eq1.0}), we have
applied in \cite{kz??}, Friedrich's extension theory \cite[pg.
28, 32]{cazS}, \cite[pg. 126-135]{zei85}: Assuming conditions
$(\mathcal{H_{\alpha}})$ or $(\mathcal{H^{\infty}_{\beta}})$, the
operator $\mathbf{T}=-\mathrm{div}(\sigma (x)\nabla\phi)$ is
positive and self adjoint, with domain of definition
$$\mathrm{D}(\mathbf{T})=\left\{\phi\in
D^{1,2}_0(\Omega,\sigma),\,\mathbf{T}\phi\in
L^2(\Omega )\right\}.$$
The space $\mathrm{D}(\mathbf{T})$, is a
Hilbert space endowed with the usual graph scalar product.  Moreover, there exist a complete system of
eigensolutions $\{e_j,\lambda_j\}$,
\begin{equation}
\left\{
\begin{array}{lll} \label{eigen1}
-\mathrm{div}(\sigma(x)\nabla e_j) = {\lambda_j} e_j, &j=1,2,..., & j\geq 1, \\
0 < {\lambda}_{1} \leq {\lambda}_{2}\leq ..., &
{\lambda}_{j}\rightarrow\infty, & \mbox{as} \; j \rightarrow
\infty.
\end{array}
\right.
\end{equation}
The fractional powers are defined as follows: For every \ $s>0, \;
\mathbf{T}^{s}$ is an unbounded selfadjoint operator in
$L^{2}(\Omega )$, with domain\
$\mathrm{D}(\mathbf{T}^{s})$ \ to be a dense subset in \ $
L^{2}(\Omega)$. The operator\ $\mathbf{T}^{s}$\ is
strictly positive and injective. Also,
$\mathrm{D}(\mathbf{T}^{s})$\ endowed with the scalar product\
$(\phi,\psi)_{\mathrm{D}(\mathbf{T}^{s})}=(\mathbf{T}^{s}\phi,\mathbf{T}^{s}\psi)_{L^{2}}$,
becomes a Hilbert space. We write as usual, \
$\mathrm{V}_{2s}=\mathrm{D}(\mathbf{T}^{s})$\ and we have the
following identifications\
$\mathrm{D}(\mathbf{T}^{-1/2})=D^{-1}_0(\Omega,\sigma)$=the dual of ${D}^{1,2}_0(\Omega,\sigma)$,
$\mathrm{D}(\mathbf{T}^{0})=L^{2}(\Omega)$\ and\
$\mathrm{D}(\mathbf{T}^{1/2})={D}^{1,2}_0(\Omega,\sigma)$.\
Moreover, the injection
$\mathrm{V}_{2s_1} \subset
\mathrm{V}_{2s_2},\;\;s_1, \, s_2 \in \mathbb{R}, \; s_1
> s_2$, is compact and dense.

While in \cite{kz??}, the local in time solvability was discussed
via compactness methods, for the purposes of the present work, it
is  more convenient to study the local in time solvability of
(\ref{eq1.0}) in ${D}^{1,2}_0(\Omega,\sigma)$, via the semigroup
method approach: The discussion above clearly shows, that the
operator $-\mathbf{T}$ is the generator of a linear strongly
continuous semigroup $\mathcal{T}(t)$ (\cite{Ball1b,cazh, Pazy83}).
\setcounter{definition}{4}
\begin{definition}
\label{defmild}
For a given function $\phi_0\in {D}^{1,2}_0(\Omega,\sigma) $, $0<\gamma <\infty$ and $T>0$, a
solution for the problem (\ref{eq1.0}), is a function $$\phi(x,t)\in
\mathrm{C}([0,T);{D}^{1,2}_0(\Omega,\sigma))\cap
\mathrm{C}^{1}([0,T);L^{2}(\Omega)),$$
satisfying
the variation of constants formula
\begin{eqnarray}
\label{pmild}
\phi(t)=\mathcal{T}(t)\phi_0+\int_{0}^{t}\mathcal{T}(t-s)f(\phi(s))
ds
\end{eqnarray}
where $f(s)=\lambda s-|s|^{2\gamma}s$.
\end{definition}
Solutions of (\ref{eq1.0}) satisfying Definition
\ref{defmild} and solutions satisfying \cite[Definition 2.3]{kz??}
(weak solutions) are the same. This is an immediate consequence of
\cite{Ball1}. \vspace{0.2cm} \\
We conclude this introductory section, by stating for the
convenience of the reader, some basic definitions and results for
our analysis. We state first a result on the existence of {\em a
branch of solutions} of an operator equation (bifurcation
in the sense  of Rabinowitz \cite{ra71}-see also \cite{dh96}).
\setcounter{theorem}{5}
\begin{theorem} \label{def1.1}
Let\ $X$\ be a Banach space with norm\ $||\cdot||_X$ and consider the operators $$\mathbf{G}(\lambda,\cdot), \mathbf{L}, \mathbf{H}(\lambda,\cdot):X\rightarrow X^*,$$ where
$\mathbf{G}(\lambda,\cdot)=\lambda \mathbf{L}(\cdot ) + \mathbf{H}(\lambda,\cdot)$, $\mathbf{L}$\
is a compact linear operator  and\ $\mathbf{H}(\lambda,\cdot)$\ is
compact and satisfies
\[
\lim_{||u||_{X} \to 0} \frac{||\mathbf{H}(\lambda,u)||_{X^*}}{||u||_{X}}=0.
\]
If\ $\lambda$\ is a simple eigenvalue of\ $\mathbf{L}$\ then the closure of
the set
\[
C=\{ (\lambda ,u) \in \mathbb{R} \times X : (\lambda,u)\;\;
\mbox{solves}\;\; \mathbf{N}(\lambda, u):=u-\mathbf{G}(\lambda,u)=0\;\;\mbox{in}\;\;X^*,\;\; u \not\equiv 0 \},
\]
possesses a maximal continuum (i.e. connected branch) of
solutions,\ $C_\lambda$,\ such that\ $(\lambda,0) \in
C_{\lambda}$\ and\ $C_{\lambda}$\ either:

(i)\ meets infinity in\ $\mathbb{R} \times X$\ or,

(ii)\ meets\ $(\lambda^*,0)$,\ where\ $\lambda^* \ne \lambda$\ is
also an eigenvalue of\ $\mathbf{L}$.
\end{theorem}

In the approximation procedure, we are making use of a generalized
Harnack-type  inequality (see \cite{dkn97, giltru77} and the
references therein).
\begin{theorem} \label{Harn} (Harnack-type Inequality) Consider
the equation
\begin{equation} \label{harn1}
- \mathrm{div}(a(x,u)\, |\nabla u|^{p-2} \nabla u ) =
f(x,u),\;\;\; x \in \Omega,
\end{equation}
where\ $\Omega \subseteq \mathbb{R}^N$,\ $1<p<N$\ and the
functions\ $a$\ and\ $f$\ satisfy the following conditions: \\
(i)\ \ $a$\ is a Carath\'{e}odory function, such that\ $a(x,s)$\
is uniformly separated from zero and bounded for almost every\ $x
\in \Omega$\ and all\ $s \in \mathbb{R}$,\  \\
(ii)\ \ $f$\ is a Carath\'{e}odory function and for any\ $M>0$\
there exists a constant\ $c_M >0$,\ such that
\[
|f(x,s)| \leq c_m |s|^{p^*-1},
\]
for almost every\ \ $x \in \Omega$\ and all\ $s \in (-M,M)$,\
where\ $p^*$\ is the critical Sobolev exponent\ $p^* =
\frac{Np}{N-p}$.

Assume that\ $u \in D^{1,p}(\Omega):=\left\{u\in L^{p^*}(\Omega)\,:\,\nabla u\in (L^p(\Omega))^N\right\}$ is a weak solution of
(\ref{harn1}) satisfying the weak formula
\[
\int_{\Omega} a(x,u) |\nabla u|^{p-2} \nabla u \nabla \phi\, dx =
\int_{\Omega} f(x,u)\, \phi\, dx,
\]
holds for any\ $\phi \in C^{\infty}_{0}(\Omega)$. Then, for any cube\ $K=K(3\rho) \subset \Omega$\
with\ $0 \leq u < M$\ in\ $K$,\ we have that
\[
\max_{x \in K_{\rho}} u(x) \leq C\; \min_{x \in K_{\rho}} u(x).
\]
In particular, if the weak solution\ $u \not\equiv 0$\ of
(\ref{eq1.1}) satisfies\ $u \geq 0$\ in\ $\Omega$\ then it follows
that\ $u$\ is strictly positive in\ $\Omega$.
\end{theorem}

\setcounter{remark}{7}
\begin{remark}
In the case where $a(x,s)\equiv a(x)$ satisfies condition (i) of Theorem \ref{harn1}, the norms of
$D_{0}^{1,2}(\Omega,a)$\ and\ $D^{1,2}(\Omega)$ are equivalent.
\end{remark}

We also recall some basic definitions and results on semiflows (see \cite{BallNS,Ball2} and \cite{jhale88,RTem88}). Let $X$ be a complete metric space. For each $\phi_0\in X$, via the correspondence $\mathcal{S}(t)\phi_0=\phi(t)$, a semiflow is a family of continuous maps $\mathcal{S}(t):X\rightarrow X$, $t\geq 0$, satisfying the semigroup identities (a) $\mathcal{S}(0)=I$, (b) $\mathcal{S}(s+t)=\mathcal{S}(s)\mathcal{S}(t)$. For $\mathcal{B}\subset X$, and $t\geq 0$
\begin{eqnarray*}
\mathcal{S}(t)\mathcal{B}:=\{\phi(t)=\mathcal{S}(t)\phi_0\;\;\mbox{with}\;\;\phi(0)=\phi_0\in\mathcal{B}\}.
\end{eqnarray*}
The {\em positive orbit} of $\phi$ through $\phi_0$ is the set $\gamma^+(\phi_0)=\{\phi(t)=\mathcal{S}(t)\phi_0,\, t\geq 0\}$. If $\mathcal{B}\subset X$ then the positive orbit of $\mathcal{B}$ is the set
\begin{eqnarray*}
\gamma^+(\mathcal{B})=\bigcup_{t\geq 0}\mathcal{S}(t)\mathcal{B}=\{\gamma^+(\phi):\,\phi(t)=\mathcal{S}(t)\phi_0\;\;\mbox{with}\;\;\phi(0)=\phi_0\in\mathcal{B}\}.
\end{eqnarray*}
If $t_0\geq 0$, $\gamma^{t_0}(\mathcal{B}):=\bigcup_{t\geq t_0}\mathcal{S}(t)\mathcal{B}=\gamma^+(\mathcal{S}(t_0)\mathcal{B})$.
The $\omega$-limit set of $\phi_0\in X$ is the set
$\omega(\phi_0)=\{z\in X\,:\,\phi(t_j)=\mathcal{S}(t_j)\phi_0\rightarrow z\;\;\mbox{for some sequence}\;\;t_j\rightarrow +\infty\}$.
A {\em complete orbit} containing $\phi_0\in X$, is a function
$\phi:\mathbb{R}\rightarrow X$ such that $\phi(0)=\phi_0$ and for
any $s\in\mathbb{R}$, $\mathcal{S}(t)\phi(s)=\phi(t+s)$ for $t\geq
0$. If $\phi$ is a complete orbit containing $\phi_0$, then the $\alpha$-limit set of
$\phi_0$ is the set $$\alpha (\phi_0)=\{z\in
X\,:\,\phi(t_j)\rightarrow z\;\;\mbox{for some
sequence}\;\;t_j\rightarrow -\infty\}.$$ The subset $\mathcal{A}$
{\em attracts} a set $\mathcal{B}$ if
$\mathrm{dist}(\mathcal{S}(t)\mathcal{B},\mathcal{A})\rightarrow
0$ as $t\rightarrow +\infty$. The set $\mathcal{A}$ is {\em
positively invariant} if
$\mathcal{S}(t)\mathcal{A}\subset\mathcal{A}$, for all $t\geq 0$
and {\em invariant} if $\mathcal{S}(t)\mathcal{A}=\mathcal{A}$ for
all $t\geq 0$. The set $\mathcal{A}$ is a global attractor if it
is compact, invariant,  and attracts all bounded sets.

The semiflow $\mathcal{S}(t)$ is {\em eventually bounded} if
given any bounded set $\mathcal{B}\subset X$, there exists
$t_0\geq 0$ such that the set $\gamma^{t_0}(\mathcal{B})$ is
bounded. The semiflow $\mathcal{S}(t)$ is said to be {\em point
dissipative} if there is a bounded set $\mathcal{B}_0$ that
attracts each point of $X$. It is called {\em asymptotically
compact} if for any bounded sequence $\phi_n$ in $X$ and for any
sequence $t_n\rightarrow\infty$, the sequence
$\mathcal{S}(t_n)\phi_n$ has a convergent subsequence. It is
called {\em asymptotically smooth} if whenever $\mathcal{B}$ is
nonempty, bounded and positively invariant, there exists a compact
set $\mathcal{K}$ which attracts $\mathcal{B}$.

A complete orbit is {\em stationary} if $\phi(t)=z$ for all $t\in\mathbb{R}$ for some $z\in X$ and each such $z$, is called an equilibrium point. We denote by $\mathcal{E}$ the set of stationary points.

The functional $\mathcal{J}:X\rightarrow\mathbb{R}$ is a {\em Lyapunov} functional for the semiflow $\mathcal{S}(t)$ if (i) $\mathcal{J}$ is continuous, (ii) $\mathcal{J}(\mathcal{S}(t)\phi_0)\leq\mathcal{J}\mathcal{S}(s)\phi_0)$ and $t\geq s\geq 0$, (iii) if $\mathcal{J}(\phi(t))=$constant for some complete orbit $\phi$ and all $t\in\mathbb{R}$, then $\phi$ is stationary.

To derive the convergence result we shall use the following Theorem.
\begin{theorem}
\label{tsibaenaball} (\cite{Ball2})\ \ Let $\mathcal{S}(t)$ be
an asymptotically compact  semiflow and suppose  that there exists a
Lyapunov function $\mathcal{J}$. Suppose further that the set
$\mathcal{E}$ is bounded. Then $\mathcal{S}(t)$ is point
dissipative, so that there exists a global attractor
$\mathcal{A}$. For each complete orbit $\phi$ containing $\phi_0$ lying in
$\mathcal{A}$ the limit sets $\alpha(\phi_0)$ and $\omega(\phi_0)$ are
connected subsets of $\mathcal{E}$ on which $\mathcal{J}$ is
constant. If $\mathcal{E}$ is totally disconnected (in particular
if $\mathcal{E}$ is countable) the limits
$$z_{-}=\lim_{t\rightarrow -\infty}\phi(t),\;\; z_{+}=\lim_{t\rightarrow+\infty}\phi (t)$$
exist and are equilibrium points. Furthermore any solution $\mathcal{S}(t)\phi_0$ tends to an equilibrium point as $t\rightarrow\infty$
\end{theorem}
\section{Global Attractor in $D^{1,2}_0(\Omega,\sigma)$}
In this section we shall show, that the degenerate semilinear
parabolic equation (\ref{eq1.0}) defines a semiflow  in the energy
space $D^{1,2}_0(\Omega,\sigma)$, possessing a global attractor.
We state first an auxiliary lemma.

\begin{lemma}
\label{boundN1}
Assume that either conditions $(\mathcal{H_{\alpha}})$
or $(\mathcal{H^{\infty}_{\beta}})$ hold. The function $f_1(s):=|s|^{2\gamma}s,\;s\in\mathbb{R}$, defines a sequentially weakly continuous map $f_1:\mathcal{D}^{1}_0(\Omega,\sigma)\rightarrow\mathrm{L^2(\Omega)}$ if
\begin{eqnarray}
\label{cruc1}
0<\gamma\leq\frac{2-\alpha}{2(N-2+\alpha)}:=\gamma^*.
\end{eqnarray}
Furthermore, if $F_1(\phi):=\int_{0}^{\phi}f_1(s)ds$, the functional $E_1:{D}^{1,2}_0(\Omega,\sigma)\rightarrow\mathbb{R}$ defined by $E_1(\phi)=\int_{\Omega}F_1(\phi)dx$,
is $C^1({D}^{1,2}_0(\Omega,\sigma),\mathbb{R})$ and  sequentially weakly continuous.
\end{lemma}
{\bf Proof:}\ \ It can be easily checked that the functional $f_1$ is well defined, under the restriction (\ref{cruc1}), by using Lemmas \ref{lemma2.1}(i)-\ref{lemma2.2}(i). Similarly, it follows that $E_1$ is well defined if
\begin{eqnarray}
\label{dontworry}
0<\gamma\leq\frac{2-\alpha}{N-2+\alpha}:=\gamma_1,
\end{eqnarray}
and note that $\gamma^*<\gamma_1$.
To show that both functionals are sequentially weakly continuous, we use the compactness of the embeddings stated in
Lemmas \ref{lemma2.1}(ii)-\ref{lemma2.2}(ii), and repeat the lines of the proof of \cite[Lemma 3.3, pg. 38 \& Theorem 3.6, pg. 40]{Ball2}.
To verify that $E_1$ is a $C^1$-functional, and its derivative is given
by
\begin{eqnarray}
\label{weak5} E_1'(\phi)(z)=\left<f_1(\phi),z\right>,\;\;\mbox{for
every}\;\;\phi\in {D}^{1,2}_0(\Omega,\sigma),\;\;z\in
{D}^{-1}_0(\Omega,\sigma),
\end{eqnarray}
we consider for $\phi,\psi\in {D}^{1,2}_0(\Omega,\sigma)$, the quantity
\begin{eqnarray}
\;\;\;\;\;\;\;\;\frac{E_1(\phi+s\psi)-E_1(\phi)}{s}&=&\frac{1}{s}
\int_{\Omega}\int_{0}^{1}\frac{d}{d\theta}F_1(\phi+\theta s\psi)d\theta dx
\nonumber\\
&=&\int_{\Omega}\int_{0}^{1}f_1(\phi+s\theta\psi)d\theta dx.
\end{eqnarray}
Setting $q=\frac{2N}{N+2-\alpha}$, $q+2_{\alpha}^*=1$, we observe that
\begin{eqnarray}
\left|\int_{\Omega}f_1(\phi+\theta s\psi)\psi dx\right|\leq c\left(\int_{\Omega}(|\phi|^{(2\gamma +1)}+|\psi|^{(2\gamma +1)})^{q}dx\right)^{\frac{1}{q}}\left(\int_{\Omega}|\psi|^{2_{\alpha}^*}dx\right)^{\frac{1}{2_{\alpha}^*}}.
\end{eqnarray}
Lemmas \ref{lemma2.1}(i)-\ref{lemma2.2}(i) are applicable under the requirement $(2\gamma +1)q\leq 2_{\alpha}^*$ which justifies (\ref{dontworry}). Using the dominated convergence theorem, we may let $s\rightarrow 0$, to obtain that $E$ is differentiable with the derivative (\ref{weak5}).

We consider next a sequence $\{\phi_n\}_{n\in\mathbb{N}}$ of ${D}^{1,2}_0(\Omega,\sigma)$ such that $\phi_n\rightarrow\phi$ in ${D}^{1,2}_0(\Omega,\sigma)$ as $n\rightarrow \infty$. It holds that
\begin{eqnarray}
\label{weak6}
\left<E_1'(\phi_n)-E_1'(\phi),z\right>\leq ||f_1(\phi_n)-f_1(\phi)||_{\mathrm{L^q}}^q||z||_{\mathrm{L^{2_{\alpha}^*}}}^{2_{\alpha}^*}.
\end{eqnarray}
Setting $p_1=\frac{2_{\alpha}^*}{q}$ we observe that the requirement for $p_1>1$, justifies the restrictions  on the exponent of degeneracy $\alpha$, imposed by $(\mathcal{H_{\alpha}})$
or $(\mathcal{H^{\infty}_{\beta}})$. Setting now $p_2=\frac{N+2-\alpha}{2(2-\alpha)}$, $p_2^{-1}+p_1^{-1}=1$, we get
\begin{eqnarray*}
||f_1(\phi_n)-f_1(\phi)||_{\mathrm{L^q}}^q&\leq& c
\left(\int_{\Omega}(|\phi_n|^{2\gamma}+|\phi|^{2\gamma})^{qp_2} dx\right)^{\frac{1}{p_2}}\left(\int_{\Omega}|\phi_m-\phi|^{2_{\alpha^*}}dx\right)^{\frac{1}{p_1}}\nonumber\\
&&:=\Lambda(\phi_m,\phi).
\end{eqnarray*}
Let $p_3=2\gamma qp_2$.  To apply Lemmas \ref{lemma2.1}(i)-\ref{lemma2.2}(i) once again, we need $p_3\leq 2_{\alpha}^*$ or (\ref{dontworry}).

Under this condition we have that $\Lambda(\phi_n,\phi)\rightarrow 0$ as $n\rightarrow\infty$ and from (\ref{weak6}), we get the continuity of $E'$.\ \ $\diamond$
\vspace{0.2cm}
\newline
We consider the energy functional $\mathcal{J}:D^{1,2}_0(\Omega,\sigma)\rightarrow\mathbb{R}$
\begin{eqnarray}
\label{Liapunov}
\mathcal{J}(\phi):=\frac{1}{2}\int_{\Omega}\sigma(x)|\nabla\phi|^2dx-\frac{\lambda}{2}\int_{\Omega}|\phi|^2dx+\frac{1}{2\gamma+2}\int_{\Omega}|\phi|^{2\gamma +2}dx.
\end{eqnarray}
\setcounter{proposition}{1}
\begin{proposition}
\label{locds} Let\ $\phi_0 \in D^{1,2}_0(\Omega,\sigma)$\ and
either conditions\ $(\mathcal{H_{\alpha}})$\ or\
$(\mathcal{H^{\infty}_{\beta}})$\ be fulfilled, and assume that (\ref{cruc1}) holds.
Then equation (\ref{eq1.0}), has a unique,
global in time (weak) solution $\phi$, such that
\begin{eqnarray}
\label{property1}
\phi\in\mathrm{C}([0,\infty);{D}^{1,2}_0(\Omega,\sigma))\cap
\mathrm{C}^1([0,\infty);L^2(\Omega)).
\end{eqnarray}
For each (weak solution) $\mathcal{J}(\phi (\cdot))\in \mathrm{C}^1([0,\infty))$ and
\begin{eqnarray}
\label{liap1}
\frac{d}{dt}\mathcal{J}(\phi(t))=-\int_{\Omega}|\partial_t\phi|^2dx
\end{eqnarray}
\end{proposition}
{\bf Proof:}\ \ By using similar arguments to those used for the
proof of Lemma \ref{boundN1}, we may show under
the assumption (\ref{cruc1}),  that the function $f(s)=|s|^{2\gamma}s-\lambda s$, defines a locally Lipschitz map $f: D^{1,2}_0(\Omega
;\sigma)\rightarrow L^2(\Omega)$.
This suffices in order to show the existence of a unique solution
$\phi$ with $\phi(0)=\phi_0$, defined on a maximal interval $[0,
T_{max})$, where $0<T_{max}\leq\infty$ \cite{cazh}.

We proceed by  showing that $T_{max}=\infty$.  First note, that by Lemma \ref{boundN1}, the energy functional $\mathcal{J}$ is $C^1$. This fact allows to adapt the method of \cite{Ball1a,Ball2}, in order to justify (\ref{liap1}) for any $t\in [0, T]$, $T<T_{max}$. We repeat the main lines of the proof, only for the shake of completeness:  For all $\phi\in\mathrm{D}(\mathbf{T})$, then
\begin{eqnarray}
\label{preliap}
\left<-\mathbf{T}\phi+f(\phi),  \mathcal{J}'(\phi)\right>&=&-\int_{\Omega}|\mathrm{div}(\sigma(x)\nabla\phi)+f(\phi)|^2dx\nonumber\\
&=&-\int_{\Omega}|\partial_t\phi|^2dx\leq 0.
\end{eqnarray}
Setting $g(t)=f(\phi(t))$ we consider  sequences $g_n(t)\in
C^{1}([0,T];D^{1,2}_0(\Omega ;\sigma))$ and $\phi_{0n}\in
\mathrm{D}(\mathbf{T})$ such that
\begin{eqnarray*}
g_n&\rightarrow& g,\;\;\mbox{in}\;\;C^{1}([0,T];D^{1,2}_0(\Omega
;\sigma)),\\
\phi_{0n}&\rightarrow&\phi_0,\;\;\mbox{in}\;\; D^{1,2}_0(\Omega
;\sigma).
\end{eqnarray*}
We define $\phi_n(t)=\mathcal{T}(t)\phi_{0n}+\int_{0}^{t}\mathcal{T}(t-s)g_n(s)ds$, and it follows from
\cite[Corrolary 2.5, p107]{Pazy83}
that $\phi_n(t)\in \mathrm{D}(\mathbf{T})$, $\phi_n\in C^{1}([0,T];D^{1,2}_0(\Omega
;\sigma))$ satisfying
$\frac{d}{dt}\phi_n(t)+\mathbf{T}\phi+f(\phi)=0$.
Also, from \cite[Lemma 5.5, pg. 246-247]{Ball1a} (see also \cite[Theorem 3.6, pg. 41]{Ball2}) we get that
\begin{eqnarray*}
\phi_n\rightarrow\phi,\;\;\mbox{in}\;\;C([0,T];D^{1,2}_0(\Omega
;\sigma)).
\end{eqnarray*}
Now using the fact that $\mathcal{J}$ is $C^1$ and (\ref{preliap}),  we may pass to the limit to
\begin{eqnarray*}
\mathcal{J}(\phi_n(t))-\mathcal{J}(\phi_{0n})&=&\int_{0}^{t}\left<\mathcal{J}'(\phi_n(s)),-\mathbf{T}\phi_n(s)+g_n(s)\right>ds\\
&=&-\int_{0}^{t}||\partial_t\phi_n(s)||^2_{\mathrm{L^2}}ds+\int_{0}^{t}\left<\mathcal{J}'(\phi_n(s)),g_n(s)-f(\phi_n(s))\right>ds
\end{eqnarray*}
to derive (\ref{liap1}).

Multiplying (\ref{eq1.0}) by ${\phi}$, and integrating over
$\Omega$, we
obtain the equation
\begin{eqnarray}
\label{eneg2} \frac{1}{2}\frac{d}{dt}||\phi||^2_{L^2}+
\int_{\Omega}\sigma(x)|\nabla\phi|^2dx-\lambda||\phi||^2_{L^2}
+\int_{\Omega}|\phi|^{2\gamma +2}dx=0.
\end{eqnarray}
We are focusing on the case where $\lambda>\lambda_1$ and the domain is unbounded.  By interpolation and
Lemma \ref{lemma2.2} (i), we have that for some $\theta\in (0,1)$,
\begin{eqnarray}
\label{unbound3}
2\lambda||\phi||^2_{L^2}&\leq&2\lambda||\phi||^{2\theta}_{L^{2\gamma+2}}
||\phi||^{2(1-\theta)}_{L^{2^*_\beta}}\nonumber\\
&\leq&2\lambda
C_{\beta}^{2(1-\theta)}||\phi||^{2\theta}_{L^{2\gamma+2}}
||\phi||^{2(1-\theta)}_{D^{1,2}_0(\Omega,\sigma)}\nonumber\\
&\leq&
\frac{1}{2}||\phi||^2_{D^{1,2}_0(\Omega,\sigma)}+c_1||\phi||^2_{L^{2\gamma+2}}\nonumber\\
&\leq& \frac{1}{2}||\phi||^2_{D^{1,2}_0(\Omega,\sigma)}+\frac{1}{2}||\phi||^{2\gamma
+2} _{L^{2\gamma+2}}+R_1,
\end{eqnarray}
where $C_\beta$ is the constant of the embedding
${D}^{1,2}_0(\Omega,\sigma)\subset
L^{2^*_\beta}(\Omega)$ and
\begin{eqnarray}
\label{numthe}
\theta&=&\frac{(\gamma
+1)(\beta-2)}{(N+\beta -2)(\gamma
+1)-N},\;\;c_1=C_\beta^2(2\lambda)^{\frac{1}{1-\theta}}(1-\theta)(2\theta)^{\frac{\theta}{1-\theta}},\nonumber\\
R_1&=&c_1^{\frac{\gamma
+1}{\gamma}}\frac{2^{\frac{1}{\gamma}}\gamma}{(\gamma +1)^{\gamma +1}}.
\end{eqnarray}
By inserting the estimate (\ref{unbound3}) to (\ref{eneg2}), we
get
\begin{eqnarray*}
\label{eneg3} \frac{1}{2}\frac{d}{dt}||\phi||^2_{L^2}+
\frac{1}{2}||\phi||_{D^{1,2}_0(\Omega,\sigma)}^2+\lambda||\phi||^2_{L^2}
+\frac{1}{2}||\phi||^{2\gamma +2}_{L^{2\gamma +2}} dx\leq
\mathrm{R}_1 .
\end{eqnarray*}
Gronwall's Lemma leads to the following inequality
 \begin{eqnarray}
 \label{abset1}
 ||\phi(t)||^2_{L^2}\leq||\phi(0)||^2_{L^2}\exp(-2\lambda t)+
\frac{1}{\lambda\mathrm{R}_1}(1-\exp(-2\lambda t)).
 \end{eqnarray}
Letting $t\rightarrow \infty$, from  (\ref{abset1}) we obtain that
\begin{eqnarray}
\label{abset2}
\limsup_{t\rightarrow\infty}||\phi(t)||^2_{L^2}\leq\rho^2,\;\;\;\rho^2=1/\lambda\mathrm{R}_1.
\end{eqnarray}
Now assume that $\phi_0$ is in a bounded set
$\mathcal{B}$ of $D^{1,2}_0(\Omega,\sigma)$. Then  (\ref{abset2}) implies that for any $\rho_1>\rho$, there exists $t_0(\mathcal{B},\rho_1)$, such that
\begin{eqnarray}
\label{boundV}
||\phi(t)||_{\mathrm{L^2}}\leq \rho_1,\;\;\mbox{for any}\;\;t\geq t_0(\mathcal{B},\rho_1).
\end{eqnarray}
We observe that by the definition of the energy functional and (\ref{boundV}),
\begin{eqnarray}
\label{abset3}
\mathcal{J}(\phi(t))&\geq& \frac{1}{2}\int_{\Omega}\sigma(x)|\nabla\phi|^2dx-\frac{\lambda}{2}\int_{\Omega}|\phi|^2dx\nonumber\\
&\geq& \frac{1}{2}\int_{\Omega}\sigma(x)|\nabla\phi|^2dx-\frac{\lambda}{2}\rho_1^2,\;\;t\geq t_0.
\end{eqnarray}
Hence, since $\mathcal{J}(\phi(t))$ is nonincreasing in $t$, we conclude that
\begin{eqnarray}
\label{boundV1}
||\phi(t)||_{D^{1,2}_0(\Omega,\sigma)}^2\leq 2\mathcal{J}(\phi_0)+\lambda\rho_1^2,\;\;t\geq t_0.
\end{eqnarray}
Thus solutions are globally defined in
$D^{1,2}_0(\Omega,\sigma)$.\ \ $\diamond$ \vspace{0.2cm}

{\bf Proof of Theorem \ref{main1}}: It is not a loss of generality to assume that $\phi_0\in\mathcal{B}(0, R)$, a closed ball of $D^{1,2}_0(\Omega,\sigma)$, of center $0$ and radius $R$. Then from Lemma \ref{cruc1} and (\ref{Liapunov}) it follows that there exists a constant $c(R)$ such that
$\mathcal{J}(\phi_0)\leq c(R)$. Hence, (\ref{boundV1}) implies  that $\mathcal{S}(t)$ is eventually bounded. Since the resolvent of the operator $-\mathbf{T}$ is compact, $\mathcal{S}(t)$ is completely continuous for $t>0$, thus asymptotically smooth. The equivalence criterion \cite[Proposition 2.3, pg. 36]{Ball2}, implies that $\mathcal{S}(t)$ is asymptotically compact. The positive orbit $\gamma^{+}(\phi_0)$ is precompact, having a nonempty compact connected invariant $\omega$-limit set $\omega (\phi_0)$.  From (\ref{liap1}) and the continuity of $\mathcal{S}(t)$ it follows that $\omega (\phi_0)\in\mathcal{E}$.

It remains to show that
$\mathcal{E}$ is bounded, to conclude that $\mathcal{S}(t)$ is point dissipative. An equilibrium point of $\mathcal{S}(t)$, is an extreme value of the functional $\mathcal{J}$ or equivalently, satisfies the weak formula
\begin{eqnarray}
\label{forstat} \int_{\Omega}\sigma(x)\nabla u\nabla v
dx-\lambda\int_{\Omega}uv dx+\int_{\Omega}|u|^{2\gamma}uv
dx,\;\;\mbox{for every}\;\;v\in D^{1,2}_0(\Omega,\sigma).
\end{eqnarray}
Setting $v=u$ in (\ref{forstat}) and using inequality (\ref{unbound3},) we obtain
\begin{eqnarray}
\label{usetocl}
\int_{\Omega}\sigma (x)|\nabla u|^2dx\leq (2\lambda)^{\frac{\gamma +1}{\gamma (1-\theta)}}R_2(\gamma ,\theta),
\end{eqnarray}
(see \ref{numthe})), which implies that for fixed $\lambda$ the set $\mathcal{E}$ is bounded.\ $\diamond$

\section{Global Bifurcation of Stationary States}
\setcounter{equation}{0} The validity of the continuous imbedding\
$D_{0}^{1,2}(\Omega,\sigma) \hookrightarrow
L^{2^{*}_{\alpha}}(\Omega)$ (Lemmas \ref{lemma2.1}-\ref{lemma2.2} (i)), enables us to use the same
arguments as in the proof \cite[Theorem 4.1 (Step 2)]{dkn97}
(see also \cite[Lemma 2.8]{abd03}), in order to prove
\ $L^{\infty}$-estimates, for the weak solutions of
(\ref{eq1.1}) and (\ref{eq1.2}).
\begin{lemma} \label{lemma2.4a}
Assume that\ $\Omega$\ is an arbitrary domain (bounded or
unbounded) and the conditions $(\mathcal{H_{\alpha}})$ or
$(\mathcal{H^{\infty}_{\beta}})$ are satisfied. Then any weak
solution\ $u$\ of (\ref{eq1.1}) or (\ref{eq1.2}) is uniformly
bounded in\ $\Omega$,\ i.e.
$||u||_{L^{\infty}(\Omega)} < C$,
where\ $C$\ is a positive constant depending on\ $\lambda$,\
$\gamma$\ and\ $K$,\ where\ $K$, is the constant appearing in
(\ref{CKNg}).
\end{lemma}

{\bf A. The bounded domain case:} We assume that the diffusion
coefficient is given by (\ref{siga}), and we consider the
following problems:
\[
(P)\;\;\;\;\;\; \biggl\{
\begin{array}{ll}
-\mathrm{div}(\sigma(x) \nabla u) = \lambda\, u
-|u|^{2 \gamma}\, u,\;\;\;\mbox{in}\;\;\Omega , \vspace{0.2cm} \\
\;\;\;\;\;\;\;\;\;\;\;\;\;\;\;\; u |_{\partial \Omega} = 0,
\end{array}
\]
\[
(P)_r\;\;\;\;\;\; \biggl\{
\begin{array}{ll}
-\mathrm{div}(\sigma(x)\nabla u) = \lambda\, u - |u|^{2 \gamma}\,
u,\;\;\;\mbox{in}\;\;\Omega_r=\Omega \backslash B_r(0),
\vspace{0.2cm} \\ \;\;\;\;\;\;\;\;\;\;\;\;\;\; u |_{\partial
\Omega_r} = 0,
\end{array}
\]
for some $r>0$ sufficiently small. Standard regularity results
(cf. \cite[Theorem 8.22]{giltru77}) imply that if\ $u$\ is a weak
solution of the problem\ $(P)$, ($(P)_r$) then\ $u \in
C^{1,\zeta}_{loc} (\Omega \backslash \{0\})$,\ ($u \in
C^{1,\zeta}_{loc} (\Omega_r)$), for some\ $ \zeta \in (0,1)$.\
\vspace{0.1cm}

For the linear eigenvalue
problems
\[
(PL),\; ((PL)_r)\;\;\;\;\;\; \biggl\{
\begin{array}{ll}
&-\mathrm{div}(\sigma(x) \nabla u) = \lambda\, u,\;\;\;\mbox{in}\;\;\Omega\;\;(\Omega_r), \vspace{0.2cm} \\
 &\;\;\;u |_{\partial \Omega} = 0,\;\;\; (u
|_{\partial \Omega_r} = 0),
\end{array}
\]
we have the following lemma.
\begin{lemma} \label{lemma2.3}
Assume that $\sigma$ is given by $(\ref{siga})$.
Problem\ $(PL)\; ((PL)_r)$, admits a positive
principal eigenvalue\ $\lambda_1$ ($\lambda_{1,r}$), given by
\begin{equation}\label{eq2.1}
\lambda_1\; (\lambda_{1,r}) = \inf_{
\begin{array}{c}
               \phi \in D_{0}^{1,2}(\Omega\; (\Omega_r),\sigma) \\
                \phi \not\equiv 0
               \end{array}}
\frac{\int_{\Omega\; (\Omega_r)} \sigma(x)\; |\nabla \phi|^2\;
dx}{\int_{\Omega\; (\Omega_r)} |\phi|^2\; dx}.
\end{equation}
with the following properties: (i) $\lambda_1$ ($\lambda_{1,r})$, is simple with a nonnegative (positive)
associated eigenfunction\ $u_1$, ($u_{1,r}$).
(ii) $\lambda_1$ ($\lambda_{1,r})$,
is the only eigenvalue of\ $(PL)\; ((PL)_r)$, with nonnegative (positive)
associated eigenfunction.
\end{lemma}
{\bf Proof:}\ The existence of $\lambda_1$ ($\lambda_{1,r}$) is a
consequence of Lemma \ref{lemma2.1} (ii) (see also
(\ref{eigen1})). For the proof of (i), let us assume that\
$u_1\geq 0$ ($u_{1,r}>0$) in $\Omega$ ($\Omega_r$) (since if\ $u$
($u_{1,r}$) is a minimizer of (\ref{eq2.1}), then\ $|u_1|$
($|u_{1,r}|$) must be also a minimizer-similar arguments may
also be find in \cite[Theorem 8.38]{giltru77}). The simplicity of
$\lambda_1$ ($\lambda_{1,r}$), can be shown by an alternative
argument, based on the so called Picone's Identity \cite{all01,
ah98, ah99}.
\newline
\vspace{0.2cm}
{$(\mathcal{PI})$}: Assume that\ $u \geq 0,\; v > 0$\ are
almost everywhere differentiable functions in\ $\Omega$. Define
\begin{eqnarray*}
L(u,v) &:=& |\nabla u|^2 + \frac{u^2}{v^2} |\nabla v|^2 -2\frac{u}{v}
\nabla u \cdot \nabla v, \\
R(u,v) &:=& |\nabla u|^2 - \nabla \biggl ( \frac{u^2}{v} \biggr )
\cdot \nabla v.
\end{eqnarray*}
Then\ $L(u,v) = R(u,v)$,\ $L(u,v) \geq 0$, and\ $L(u,v)= 0$, if
and only if\ $u=kv$\ for some constant\ $k$,\ a.e. in\ $\Omega$.\
\vspace{0.1cm}

Let $\Omega_0 \subset \Omega$ a compact subset of $\Omega$,  and
$0 \leq \phi \in C_{0}^{\infty}(\Omega)$.  For $\lambda>0$ we
consider $u \in C^{1,\zeta}_{loc}(\Omega)$, $ \zeta \in (0,1)$, a
weak solution of $(PL)$, such that $0\leq u$ a.e in $\Omega$.
Then, for any\ $\varepsilon>0$,\ we have that
\begin{eqnarray} \label{eq2.1a}
0 &\leq& \int_{\Omega_0} \sigma(x)\; L(\phi, u+\epsilon)\; dx \leq
\int_{\Omega} \sigma(x)\; L(\phi, u+\epsilon)\; dx = \nonumber \\
&=& \int_{\Omega} \sigma(x)\; R(\phi, u+\epsilon)\; dx = \nonumber \\
&=& \int_{\Omega} \sigma(x)\; |\nabla \phi|^2\; dx - \int_{\Omega}
\sigma(x)\; \nabla \biggl ( \frac{\phi^2}{u+\epsilon} \biggr )
\cdot \nabla
u\; dx = \nonumber \\
&=& \int_{\Omega} \sigma(x)\; |\nabla \phi|^2\; dx + \int_{\Omega}
\biggl ( \frac{\phi^2}{u+\epsilon} \biggr )  \nabla(\sigma(x)
\nabla u)\; dx = \nonumber \\
&=& \int_{\Omega} \sigma(x)\; |\nabla \phi|^2\; dx - \lambda
\int_{\Omega} \biggl ( \frac{\phi^2}{u+\epsilon} \biggr ) u\; dx.
\end{eqnarray}

Assume now that\ $\lambda_1$\ is not simple. Let\ $v \not\equiv
u_1$\ be another associated eigenfunction,\ $v \in
D_{0}^{1,2}(\Omega,\sigma)$\ almost everywhere differentiable in\
$\Omega$,\ such that\ $v(x) \geq 0$\ in some\ $\Omega^+ \subset
\Omega$. Consider (\ref{eq2.1a}) with\ $\Omega_0 \subseteq
\Omega^+$,\ $\lambda=\lambda_1$\ and\ $u=u_1$.\ Letting\ $\phi \to
v$\ in\ $\Omega^+$\ and\ $\epsilon \to 0$,\ Fatou's Lemma and
Lebesgue Dominated Convergence Theorem, imply that\ $L(v,u_1) = 0$
a.e. in\ $\Omega^+$. Hence from $(\mathcal{PI})$ we get that\ $v=k
u_1$,\ a.e. in\ $\Omega^+$, which implies the simplicity of\
$\lambda_1$. Property (i) is proved.

For the proof of (ii), we suppose that there exists another
eigenvalue of $(PL)$, $\lambda^* >\lambda_1$,\ to which
corresponds a nonnegative eigenfunction\ $u^*$.\ Consider
(\ref{eq2.1a}) with\ $\Omega_0 \subseteq \Omega^+$,\
$\lambda=\lambda^*$\ and\ $u=u_*$.\ Letting\ $\phi \to u_1$\ in\
$\Omega$\ and\ $\epsilon \to 0$,\ we obtain that
\[
0 \leq \int_{\Omega} \sigma(x)\; L(u_1, u^*)\; dx <0,
\]
which is a contradiction.\ $\diamond$
\begin{lemma} \label{lemma2.4}
Assume that $\sigma$  is given by (\ref{siga}). Let also \
$\lambda_1$, $\lambda_{1,r}$, be the positive principal
eigenvalues of the problems $(PL)$, $(PL)_r$, respectively. \
Then,\ $u_{1,r} \to u_1$\ in\ $D_{0}^{1,2}(\Omega, \sigma) \cap
L^{\infty}_{loc} (\Omega\setminus\{0\})$, and\ $\lambda_{1,r}
\downarrow \lambda_1$,\ as\ $r \downarrow 0$.\
\end{lemma}
{\bf Proof:}\ We extend  $u_{1,r}$ on $\Omega$ as
\[
\hat{u}_{1,r} (x) = : \left\{
\begin{array}{ll}
u_{1,r} (x),\;\;\;& x \in \Omega_r, \\
0,\;\;\;& x \in B_r,
\end{array}
\right.
\]
for any sufficiently small\ $r>0$, but in the sequel, for
convenience, we shall use the same notation
$u_{1,r}\equiv\hat{u}_{1,r}$.\ Observe that
\[
\lambda_{1,r} = \frac{\int_{\Omega_r} \sigma(x)\; |\nabla
u_{1,r}|^2\; dx}{\int_{\Omega_r} |u_{1,r}|^2\; dx} =
\frac{\int_{\Omega} \sigma(x)\; |\nabla u_{1,r}|^2\;
dx}{\int_{\Omega} |u_{1,r}|^2\; dx} \geq \lambda_1,
\]
and\ $\lambda_{1,r}$\ is an decreasing sequence, as\ $r \to 0$,\
since\ $\Omega_{\rho} \subset \Omega_{\varrho}$,\ for any\ $\rho >
\varrho$ .\ Clearly, $u_{1,r}$ forms  a bounded sequence in\
$D_{0}^{1,2}(\Omega, \sigma)$.  Lemma \ref{lemma2.1} (ii) and
Lemma \ref{lemma2.3}, imply the existence of a pair\
$(\lambda^*,u^*)$, and a subsequence of $u_{1,r}$ (not
relabelled), such that
\[
\lambda_{1,r} \int_{\Omega} |u_{1,r}|^2dx \to \lambda^*
\int_{\Omega} |u^*|^2dx,
\]
as\ $r \to 0$.\ Then (\ref{eq2.1}) implies  that $u_{1,r} \to
u^*$\ in\ $D_{0}^{1,2}(\Omega, \sigma)$\ and\ $u^*$\ satisfies
\[
\int_{\Omega} \sigma(x)\; |\nabla u^*|^2\; dx = \lambda^*
\int_{\Omega} |u^*|^2\; dx.
\]
From Lemma \ref{lemma2.3} (ii), we obtain that \ $(\lambda^*,u^*)
\equiv (\lambda_1,u_1)$.  We conclude by justifying the claim
that\ $u_{1,r}$ is uniformly bounded in\ $L^{\infty}(\Omega)$.
Note that\ $\lambda_{1,r} \in (\lambda_1, \lambda_1 + \epsilon)$,\
for some\ $\epsilon >0$\ and any\ $r$\ small enough. Since\
$u_{1,r} \in D_{0}^{1,2}(\Omega, \sigma)$, it holds that\
$||u_{1,r}||_{L^{2^{*}_{\alpha}}(\Omega)} < K\,
||u_{1,r}||_{D_{0}^{1,2}(\Omega,\sigma)}$,\ $K$\ is given in
(\ref{CKNg}) and is independent of\ $r$.\ Hence, from Lemma
\ref{lemma2.4a} we have that\ $u_{1,r}$\ is uniformly bounded in\
$L^{\infty}(\Omega)$.\ Then by a standard bootstrap argument we
get that\ $u_{1,r} \to u_1$\ in\ $L^{\infty}_{loc}
(\Omega\setminus\{0\})$\ and the proof is completed.\
\vspace{0.2cm} $\diamond$

\setcounter{proposition}{3}
\begin{proposition}
\label{thm2.1} Assume that $\sigma$ is given by (\ref{siga}). The principal eigenvalues\ $\lambda_1,\;
\lambda_{1,r}$\ of the linear problems\ $(PL),\; (PL)_r$ , are
bifurcation points of the problems\
$(P),\; (P)_r$ respectively.  Moreover, for any (sufficiently small) $r>0$, the branch\ $C_{\lambda_{1,r}}$\ is
global, and any function which belongs to $C_{\lambda_{1,r}}$, is strictly positive.
\end{proposition}
{\bf Proof:}\ The existence of branches bifurcating from\
$\lambda_1,\; \lambda_{1,r}$\ follows by Theorem \ref{def1.1},
since  Lemma \ref{lemma2.1} (ii) and  Lemmas
\ref{lemma2.3}-\ref{lemma2.4}, are in hand.  
We outline the proof for the branch $C_{\lambda_1}$.

As in \cite{bro94}, we define a bilinear form  in $C^{\infty}_0(\Omega)$ by
\begin{eqnarray}
\label{Rabi1}
\left<u,v\right>=\int_{\Omega}\sigma (x)\nabla u\nabla v\,dx-\frac{c^{-1}}{2}\int_{\Omega}uv\,dx,\;\;\mbox{for all}\;\;u,v\in C^{\infty}_0(\Omega).
\end{eqnarray}
($c$ is the constant in (\ref{poi})) and we define $X$ to be the completion of $C^{\infty}_0(\Omega)$ with respect to the norm induced by  (\ref{Rabi1}),
$||u||_{X}^2=\left<u,u\right>$: from inequality \ref{poi} we get that
\begin{eqnarray*}
\frac{1}{2}||u||^2_{D^{1,2}_{0}(\Omega,\sigma)}\leq ||u||_{X}^2\leq\frac{3}{2}||u||^2_{D^{1,2}_{0}(\Omega,\sigma)},
\;\;\mbox{for all}\;\;u,v\in C^{\infty}_0(\Omega),
\end{eqnarray*} 
and by density it follows that $X=D^{1,2}_{0}(\Omega,\sigma)$. Henceforth we may suppose that the norm in $X$ coincides with the norm in $D^{1,2}_{0}(\Omega,\sigma)$ and that the inner product in $X$ is given by $<u,v>=(u,v)_{\sigma}$ (moreover, we may assume that if $<\cdot,\cdot>_{X,X^*}$ denotes the duality pairing on $X$, then $<\cdot,\cdot>_{X,X^*}=<\cdot,\cdot> $ \cite[Identification Principle 21.18, pg. 254]{zei85}). On the other hand, the bilinear form
\begin{eqnarray*}
\mathbf{a}(u,v)=\int_{\Omega}uv\,dx,\;\;\mbox{for all}\;\;u,v\in X,
\end{eqnarray*}
is clearly continuous in $X$ as it follows from Lemma \ref{lemma2.1} , and by the Riesz reperesentation theorem we can define a bounded linear operator $\mathbf{L}$ such that
\begin{eqnarray}
\label{Rabi1.1}
\mathbf{a}(u,v)=<\mathbf{L}u,v>,\;\;\mbox{for all}\;\;u,v\in X.
\end{eqnarray}
The operator $\mathbf{L}$ is self adjoint and by Lemma \ref{lemma2.1} (ii) is compact. Th largest eigenvalue $\nu_1$ of $\mathbf{L}$ is given by 
\begin{eqnarray*}
\nu_1=\sup_{u\in X}\frac{<\mathbf{L}u,u>}{<u,u>}=\sup_{u\in X}\frac{\int_{\Omega}u^2dx}{\int_{\Omega}\sigma (x)|\nabla u|^2dx}.
\end{eqnarray*}
It follows from Lemma \ref{lemma2.3} that the the positive eigenfunction $u_1$ of $(PL)$ corresponding to $\lambda_1$ is a positive eigenfunction of $\mathbf{L}$ corresponding to $\nu_1=1/\lambda_1$.
We consider now the nonlinear operator $\mathbf{N}(\lambda,\cdot):\mathbb{R}\times X\rightarrow X^*$ defined by
\begin{eqnarray}
\label{Rabi2}
\left<\mathbf{N}(\lambda,u),v\right>=\int_{\Omega}\sigma (x)\nabla u\nabla v\,dx-\lambda\int_{\Omega}uv\,dx+\int_{\Omega}|u|^{2\gamma}uv\,dx,\;\;\mbox{for all}\;\;v\in X.
\end{eqnarray}
Arguments very similar to those used for the proof of Lemma \ref{boundN1}, can be used in order to verify that for fixed $u\in X$, the functional $S$ defined by 
\begin{eqnarray*}
S(v)= \int_{\Omega}\sigma (x)\nabla u\nabla v\,dx-\lambda\int_{\Omega}uv\,dx+\int_{\Omega}|u|^{2\gamma}uv\,dx,\;v\in X,
\end{eqnarray*}
is a bounded linear functional and thus $\mathbf{N}(\lambda, u)$ is well defined from (\ref{Rabi2}). Moreover  by using the fact that $X=D^{1,2}_0(\Omega,\sigma)$ and relation (\ref{Rabi1.1}), we can rewrite $\mathbf{N}(\lambda,u)$ in the form $\mathbf{N}(\lambda,u)=u-\mathbf{G}(\lambda,u)$
where $\mathbf{G}(\lambda,u):=\lambda\mathbf{L}u-\mathbf{H}(u)$,
\begin{eqnarray*}
<\mathbf{H}(u),v>=\int_{\Omega}|u|^{2\gamma}uv\,dx\;\;\mbox{for all}\;\;v\in X.
\end{eqnarray*}
The restriction (\ref{cruc1}) and  Lemma \ref{lemma2.1} (ii) implies that $\mathbf{H}$ is compact. Moreover we observe that
\begin{eqnarray}
\label{limbif}
\frac{1}{||u||_{X}}|<\mathbf{H}(u),v> |&\leq&
\frac{1}{||u||_{X}}||u||^{2\gamma}_{L^{2\gamma+2}}
||u||_{L^{2\gamma+2}}
||v||_{L^{2\gamma+2}}\nonumber\\
&\leq& c_1\, ||u||^{2\gamma}_{X}
||v||_{X}.
\end{eqnarray}
Therefore, we get from (\ref{limbif}) that
\begin{eqnarray*}
\lim_{||u||_{X}\rightarrow 0}\frac{||\mathbf{H}(u)||_{X^*}}{||u||_{X}}
=\lim_{{||u||_{X}\rightarrow 0}}\sup_{||v||_{X}\leq1}\frac{1}{||u||_{X}}|<\mathbf{H}(u),v>|=0.
\end{eqnarray*}

To prove that\
$C_{\lambda_{1,r}}$\ is global for sufficiently small $r>0$, we
proceed in two steps. \vspace{0.2cm} \\
(a)\ We shall prove first that for all solutions\ $(\lambda,u) \in
C_{\lambda_1,r}$\ close to\ $(\lambda_{1,r} ,0)$ it holds that
$u(x)>0$,\ $x \in \Omega_r$. In other words,  we have to show that
there exists\ ${\epsilon}_0>0$,\ such that\ for any\ $(\lambda,
u(x)) \in C_{\lambda_{1,r}} \cap B_{\epsilon_0}
((\lambda_{1,r},0))$, it holds that $u(x) > 0$,\ for any\ $x \in
\Omega_r$\ (By $B_{\epsilon_0} ((\lambda_{1,r},0))$, we denote the
open ball of $C_{\lambda_{1,r}}$ of center $(\lambda_{1,r},0)$ and
radius $\epsilon_0$).

We argue by contradiction: Let\ $(\lambda_n, u_n)$ \ be a sequence
of solutions of\ $(P)_r$,\ such that \ $(\lambda_n, u_n) \to
(\lambda_{1,r},0)$\ and assume that $u_n$\ are changing sign in\
$\Omega_r$. 
Let $u_{n}^{-} := \min \{0, u_n\}$ and
${\mathcal{U}}_n^{-} =: \{ x \in \Omega_r: u_n(x)<0 \}$.\ 
Since $u_n=u_n^+ -u_n^{-}$ is a solution of the problem $(P)_r$ it can be easily seen that $u_n^{-}$, satisfies (in the weak sense) the equation
\begin{eqnarray}
\label{Prtonos}
-\mathrm{div}\left(\sigma(x)\nabla u_n^{-}\right)&-&\lambda_n u_n^{-}+|u_n|^{2\gamma}u_n^{-}=0,\\
u_n^{-}|_{\partial \Omega_r}&=&0.\nonumber
\end{eqnarray}

Then,
multiplying $(\ref{Prtonos})$ with $u_{n}^{-}$ and integrating over $\Omega_r$ we have that
\begin{eqnarray}
\label{Prtonos2}
\int_{{\mathcal U}_n^{-}} | \nabla u_{n}^{-} |^{2}\, dx  -
\lambda_n \int_{{\mathcal U}_n^{-}} |u_{n}^-|^2\, dx+\int_{{\mathcal U}_n^{-}}|u_n|^2|u_n^{-}|^2dx=0.
\end{eqnarray}
Since $\lambda_n$ is a bounded sequence, it follows  from (\ref{Prtonos2}), H\"{o}lder's inequality and relation (\ref{CKNg}) that
\begin{eqnarray*}
||u_{n}^{-}||^{2}_{D_{0}^{1,2}({\mathcal U}_n^{-},\sigma)}  &\leq& \lambda_n\int_{{\mathcal U}_n^{-}} |u_{n}^-|^2\, dx\\
&\leq& C\,|{\mathcal U}_n^{-}|^{\frac{2-\alpha}{N}}\left(\int_{{\mathcal U}_n^{-}}|u_n|^{2^*_{\alpha}}\right)^{\frac{2}{2^*_{\alpha}}}\\
&\leq&
C |{\mathcal U}_n^{-}|^{\frac{2-\alpha}{N}}||u_{n}^-||^{2}_{D_{0}^{1,2}({\mathcal U}_n^{-},\sigma)}. 
\end{eqnarray*}
or, equivalently
\begin{equation} \label{sob}
M \leq |{\mathcal U}_n^{-}|,\;\;\mbox{for all}\;\; n,
\end{equation}
where the constant \ $M$ \ is independent of \ $n$.
We denote now  by\ $\tilde{u}_n,=u_n/||u_n||$\ the
normalization of\ $u_n$.\ Then there exists a subsequence of\
$\tilde{u}_n$ (not relabelled) converging weakly in\
$D_{0}^{1,2}(\Omega_r,\sigma)$ \ to some function\ $\tilde{u}_0$.
By following the lines of the proof of Lemma \ref{lemma2.4}, it
can be seen that $\tilde{u}_0=u_{1,r}$. 
Moreover,
$\tilde{u}_n\to u_{1,r} >0$ in $L^2(\Omega_r)$. Passing to a further subsequence if necessary,  by Egorov's
Theorem, $\tilde{u}_n \to u_{1,r} $ uniformly on $\Omega_r$ with the exception of a set of arbitrary small measure.
This contradicts (\ref{sob}) and we conclude
the functions \ $u_n$ cannot change sign
(for a similar argument, we refer to
\cite{dh96, dkn97, zs02}). \vspace{0.2cm} \\
(b)\ Suppose now that for some solution\ $(\lambda,u) \in
C_{\lambda_1,r}$, there exists a point\ $x_0 \in \Omega_r$,\ such
that\ $u(x_0) < 0$. Using (a), the fact that the continuum
$C_{\lambda_1,r}$\ is connected (see Theorem \ref{def1.1}) and
the\ $C^{1,\zeta}_{loc}(\Omega_r)$- regularity of solutions, we
get that there exists\ $({\lambda}_0,u_0) \in C_{\lambda_1,r}$,\
such that\ $u_0(x) \geq 0$,\ for all\ $x \in \Omega_r$,\ except
possibly some point\ $x_0 \in \Omega_r$,\ such that\ $u_0(x_0) =
0$.\ Then Theorem \ref{Harn}, implies that\ $u_0 \equiv 0$\ on\
$\Omega_r$. Thus, we may construct a sequence\ $\{
({\lambda}_n,u_n) \} \subseteq C_{\lambda_1,r}$,\ such that\
$u_n(x)> 0$,\ for all\ $n$\ and\ $x \in \Omega_r$,\ $u_n
\rightarrow 0$\ in\ $D^{1,2}_0(\Omega_r, \sigma)$,\ and\
${\lambda}_n \rightarrow {\lambda}_0$.\ However, this is true only
for\ $\lambda_0 = \lambda_1$.\ As a consequence, we have that\
$C_{\lambda_1,r}$\ cannot cross $(\lambda,0)$ for some $\lambda
\ne \lambda_1$, and every function which belongs to $C_{\lambda_1,r}$
is strictly positive. $\diamond$
\setcounter{theorem}{4}
\begin{theorem} \label{thm2.2}
Assume that $\sigma$ is given by (\ref{siga}).
Then,\ $C_{\lambda_1}$\ is a global branch of nonnegative
solutions for the problem\ $(P)$.\
\end{theorem}
\emph{Proof}\ It suffices to prove that\ $C_{\lambda_{1,r}} \to
C_{\lambda_{1}}$,\ as\ $r \to 0$. The global character of\
$C_{\lambda_{1,r}}$\ implies that for any fixed positive number\
$R$,\ and any\ $r$\ sufficiently small, the set $C_{\lambda_{1,r}}
\cap B_R (\lambda_{1,r},0)$\ is not empty. By using the properties
of\ $\lambda_1$ established in Lemma \ref{lemma2.3}\ and the
compactness arguments of Lemma \ref{lemma2.4}, we can show that
\[
\lim_{r \to 0} C_{\lambda_{1,r}} \cap B_R (\lambda_{1,r},0) \to
C_{\lambda_{1}} \cap B_R (\lambda_{1},0),\;\;\; \mbox{for every}\;
R>0,
\]
which implies that\ $C_{\lambda_{1,r}} \to C_{\lambda_{1}}$,\ as\
$r \to 0$.\ Alternatively, one may use Whyburn's Theorem
\cite{amb97, bl04, eg00, giac98}.\ $\diamond$ \vspace{0.2cm}

{\bf Proof of Theorem \ref{main2} in the case of
$(\mathcal{H_{\alpha}})$:} One has to extend Theorem \ref{thm2.2}
in the case of a diffusion coefficient satisfying
$(\mathcal{H_{\alpha}})$. Since the set of zeroes of $\sigma$,
$Z_{\sigma}$ is finite, we  may use (\ref{charw}) and consider
approximating problems similar to $(P_r)$, defined this time in
the domain $\Omega_r=\Omega\setminus\bigcup_i B_r(z_i)$. The
finiteness of  $Z_{\sigma}$, allows to repeat the proofs of Lemmas
\ref{lemma2.3}-\ref{lemma2.4} and Proposition \ref{thm2.1},
without additional complications. \vspace{0.3cm} \\
{\bf B. The unbounded domain case} We assume that the diffusion
coefficient is given by (\ref{sigab}) and we consider the
following problem:
\[
(P)_{\infty}\;\;\;\;\;\; \biggl\{
\begin{array}{ll}
&-\mathrm{div}(\sigma(x) \nabla u) = \lambda\, u
-|u|^{2 \gamma}\, u,\;\;\mbox{in}\;\;\Omega, \vspace{0.2cm} \\
&\;\;\;\;\;\;\;\;\;\;\;\;\;\;\;\;u|_{\partial
\Omega} = 0.
\end{array}
\]
where\ $\Omega \subseteq \mathbb{R}^N$,\ $N \geq 2$,\ is an
unbounded domain containing the origin. The regularity
results of \cite[Theorem 8.22]{giltru77}, imply once again that
if\ $u$\ is a weak solution of the problem\ $(P_{\infty})$, then\
$u \in C^{1,\zeta}_{loc} (\Omega \backslash \{0\})$, for some\ $
\zeta \in (0,1)$.\ This time, we consider the approximating
problem,
\[
(P)_R\;\;\;\;\;\; \biggl\{
\begin{array}{ll}
&-\mathrm{div}(\sigma(x) \nabla u) = \lambda\, u
-|u|^{2 \gamma}u,\;\;\mbox{in}\;\;\Omega_R = \Omega \cap B_R(0), \vspace{0.2cm} \\
&\;\;\;\;\;\;\;\;\;\;\;\;\;\;u |_{\partial
\Omega_R} = 0.
\end{array}
\]
We consider the linear
eigenvalue problems
\[
(PL)_{\infty},\; ((PL)_R)\;\;\;\;\;\; \biggl\{
\begin{array}{ll}
-\mathrm{div}(\sigma(x) \nabla u) = \lambda\, u,\;\;\mbox{in}\;\;\Omega\;\;(\Omega_R),
\vspace{0.2cm} \\ \;\;\;\;\;\;\;\;\;\;\;\;\;\;\;\; u|_{\partial \Omega} =
0,\;\;\; (u |_{\partial \Omega_R} = 0).
\end{array}
\]
A result similar to Lemma \ref{lemma2.3}, holds.
\setcounter{lemma}{5}
\begin{lemma} \label{lemma2.7}
Assume that $\sigma$ is given by $(\ref{sigab})$.
Problem\ $(PL)_{\infty}\; ((PL)_R)$, admits a positive
principal eigenvalue\ $\lambda_1$ ($\lambda_{1,R}$), given by
\begin{equation}\label{eq2.2}
\lambda_1\; (\lambda_{1,R}) = \inf_{
\begin{array}{c}
               \phi \in D_{0}^{1,2}(\Omega\; (\Omega_R),\sigma) \\
                \phi \not\equiv 0
               \end{array}}
\frac{\int_{\Omega\; (\Omega_R)} \sigma(x)\; |\nabla \phi|^2\;
dx}{\int_{\Omega\; (\Omega_R)} |\phi|^2\; dx}.
\end{equation}
with the following properties: (i) $\lambda_1$ ($\lambda_{1,R})$, is simple with a nonnegative (positive)
associated eigenfunction\ $u_1$, ($u_{1,R}$).
(ii) $\lambda_1$ ($\lambda_{1,R})$,
is the only eigenvalue of\ $(PL)_{\infty}\; ((PL)_R)$, with nonnegative (positive)
associated eigenfunction.
\end{lemma}
To prove a similar to Lemma \ref{lemma2.4} result, we shall use the extension
\[
\hat{u}_{1,R} (x) = : \left\{
\begin{array}{ll}
u_{1,R} (x),\;\;\;& x \in \Omega_R, \\
0,\;\;\;& x \in \Omega \backslash \Omega_R,
\end{array}
\right.
\]
and use for convenience the notation $\hat{u}_{1,R}\equiv u_{1,R}$.
\begin{lemma} \label{lemma2.8}
Let\ $\lambda_1$, ($\lambda_{1,R}$)\ be the positive principal
eigenvalues of the problems\ $(PL_{\infty}),\; ((PL)_R)$.\ Then,\
$\hat{u}_{1,R} \to u_1$,\ in\ $D_{0}^{1,2}(\Omega, \sigma) \cap
L^{\infty}_{loc} (\Omega\setminus\{0\})$\ and\ $\lambda_{1,R} \downarrow
\lambda_1$,\ as\ $R \to \infty$.\
\end{lemma}
We remark that for each $R>0$, Theorem \ref{thm2.2} is applicable
for $(P)_R$: There exists a global branch,\ $C_{\lambda_{1,R}}$,\
of nonnegative solutions, bifurcating from\ $\lambda_{1,R}$. This
suffices for a repetition of arguments similar to those used for
the proof of Theorem \ref{thm2.2}, to show that\
$C_{\lambda_{1,R}} \to C_{\lambda_1}$,\ as\ $R \to \infty$.\
\setcounter{theorem}{7}
\begin{theorem} \label{thm2.3}
Assume that $\sigma$ is given by (\ref{sigab}). Then,\ $\lambda_1$\ is a
bifurcating point of the problem\ $(P_{\infty})$\ and\
$C_{\lambda_1}$\ is a global branch of nonnegative solutions.
\end{theorem}

{\bf Proof of Theorem \ref{main2} in the case of
$(\mathcal{H_{\beta}^{\infty}})$:} One has to consider
approximating problems similar to $(P_R)$, defined in the domain
$\Omega_R=\Omega \cap B_R(0)$. The conclusion follows from Theorem
\ref{main2} in the case of $(\mathcal{H_{\alpha}})$,\ repeating
the proofs of Lemmas \ref{lemma2.7}-\ref{lemma2.8} and the
arguments of the proof of Theorem \ref{thm2.2} .
\vspace{0.3cm} \\
{\bf C. Properties of the global branches} In the remaining part
of this section, we state some further properties of the global
branch\ $C_{\lambda_{1}}$,\ both in the bounded and the unbounded
domain case. For similar properties possessed by solutions of nondegenerate
elliptic equations, we refer to  \cite{amb97, eg00, giac98}.
\setcounter{lemma}{8}
\begin{lemma} \label{lemma2.9}
Assume that\ $\Omega$\ is a bounded domain and the condition
$(\mathcal{H_{\alpha}})$ is fulfilled, $\lambda>\lambda_1$ is
fixed and $u_{\lambda,r} \in C_{\lambda_{1,r}}$\ and\ $u_{\lambda}
\in C_{\lambda_{1}}$.\ Then, we have that
\[
u_{\lambda,r}(x) \leq u_{\lambda}(x),\;\;\; \mbox{for any}\; x \in
\bar{\Omega}_r,\;\; \mbox{and any}\; r \to 0,
\]
and
\[
u_{\lambda,r} \to u_{\lambda}\;\; \mbox{in}\;
L^{\infty}_{loc}(\Omega\setminus\{0\}),\;\; \mbox{as}\; r \to 0.
\]
\end{lemma}
{\bf Proof:}\ The solution $u_{\lambda,r}$ satisfies
$(P)_r$,\ while\ $u_{\lambda}$\ satisfies
\[
\biggl\{
\begin{array}{ll}
-\mathrm{div}(\sigma(x) \nabla u) = \lambda\, u
-|u|^{2 \gamma}\, u, \vspace{0.2cm} \\
\;\;\;\;\;\; u |_{\partial \Omega} = 0,\;\;\;\; u|_{\partial B_r}
>0.
\end{array}
\]
Having in mind, that both\ $u_{\lambda,r}$\ and\ $u_{\lambda}$\ are
sufficiently smooth and positive functions on\ $\bar{\Omega}_r$, from
the Comparison Principle \cite[Theorem 10.5]{pucc04},  we conclude
the first assertion of Lemma.

Next, we proceed as in Lemma \ref{lemma2.4}. Since\ $u_{\lambda,r}
\in D_{0}^{1,2}(\Omega, \sigma)$, it holds that
$||u_{\lambda,r}||_{L^{2^{*}_{\alpha}}(\Omega)} < K\,
||u_{\lambda,r}||_{D_{0}^{1,2}(\Omega,\sigma)}$, where $K$ is
given in (\ref{CKNg}) and is independent of $r$.  Hence, from
Lemma \ref{lemma2.4a} we have that  $u_{\lambda,r}$ is uniformly
bounded in\ $L^{\infty}(\Omega)$.\ Consider\ $\psi = u -
u_{\lambda,r}$.\ Then, from \cite[Theorem 8.8]{giltru77} we obtain
that
\[
||\psi||_{W^{2,2}_{loc} (\Omega\setminus\{0\})} \leq c\,
||\psi||_{W^{1,2} (\Omega)} + O(r),\;\;\; \mbox{as}\;\; r \to 0,
\]
for some positive constant\ $c$\ independent from\ $r$.\ Then, by
a standard bootstrap argument, we conclude that $u_{\lambda,r} \to
u_{\lambda}$ in\ $L^{\infty}_{loc} (\Omega\setminus\{0\})$\ and
the proof is completed.\ \vspace{0.2cm} $\diamond$

Similar results may be obtained for the unbounded domain case.
\setcounter{lemma}{9}
\begin{lemma}
Assume that\ $\Omega$\ is an unbounded domain and the condition
$(\mathcal{H_{\beta}^{\infty}})$ is fulfilled, $\lambda>\lambda_1$
is fixed, $u_{\lambda,R} \in C_{\lambda_{1,R}}$\ and\ $u_{\lambda}
\in C_{\lambda_{1}}$.\ Then, we have that
\[
u_{\lambda,R}(x) < u_{\lambda}(x),\;\;\; \mbox{for any}\; x \in
\Omega_R,
\]
and
\[
u_{\lambda,R} \to u_{\lambda}\;\; \mbox{in}\;
L^{\infty}_{loc}(\Omega\setminus\{0\}),\;\; \mbox{as}\; R \to \infty.
\]
\end{lemma}
For both, the bounded and the unbounded domain case, we have the
following
\setcounter{proposition}{10}
\begin{proposition}
Assume that conditions $(\mathcal{H_{\alpha}})$ or $(\mathcal{H_{\beta}^{\infty}})$ hold.  Then,

(i)\ \ \ The global branch\ $C_{\lambda_{1}}$\ bends to the right
of\ $\lambda_1$\ (supercritical bifurcation) and it is bounded
for\ $\lambda$\ bounded.

(ii)\ \ Every solution\ $u \in C_{\lambda_{1}}$,\ is the unique
nonnegative solution for the problem (\ref{eq1.1}).
\end{proposition}
{\bf Proof:}\, (i)\ Assume that\ $C_{\lambda_{1}}$\ bends to the
left of\ $\lambda_{1}$.\ Then there exists a pair\ $(\lambda,u)
\in \mathbb{R} \times D_{0}^{1,2}(\Omega,\sigma)$,\
$0<\lambda<\lambda_1$,\ such that
\begin{equation} \label{eq2.1ab}
\int_{\Omega} \sigma(x) |\nabla u|^2\; dx = \lambda \int_{\Omega}
|u|^2\; dx - \int_{\Omega} |u|^{2\gamma+2}\; dx,
\end{equation}
The last equality implies that
\[
||u||^{2}_{D_{0}^{1,2}(\Omega,\sigma)} \leq \lambda
||u||^{2}_{L^2(\Omega)},\;\;\; \mbox{with}\;\; \lambda<\lambda_1,
\]
which contradicts the variational characterization
(\ref{eq2.1}) of\ $\lambda_1$.\ Thus,\ $C_{\lambda_{1}}$\ must
bend to the right of\ $\lambda_{1}$.\ To show  that\
$C_{\lambda_{1}}$\ is bounded for\ $\lambda$\ bounded, we proceed
exactly as for the derivation of the estimate (\ref{usetocl}).

(ii)\ Let\ $u \in C_{\lambda_1}$,\ and suppose that\ $v$\ is a
nonnegative solution of (\ref{eq1.1}) with\ $u \not\equiv v$.\ We
claim that\ $u(x) \leq v(x)$,\ for any\ $x \in \Omega\setminus\{0\}$.\
This is a concequence of Lemma \ref{lemma2.9}, since by
the Comparison Principle we have that
\[
u_{\lambda,r}(x)\;\; (\mbox{or}\; u_{\lambda,R}(x))  \leq \min_{x
\in \Omega} \{u(x),\; v(x)\},
\]
and of  the\ $L^{\infty}_{loc}$-convergence of\ $u_{\lambda,r}$\
(or\ $u_{\lambda,R}(x)$) to\ $u$.\ Then, from (\ref{eq2.1ab}) we
must have that
\[
\int_{\Omega} (|u|^{2\gamma}v -
|v|^{2\gamma}u)\, dx = 0,
\]
which is a contradiction, unless\ $u \equiv v$.\ $\diamond$

We  emphasize  that  uniqueness results in the
case of semilinear elliptic equations, have been treated by many authors.
We refer to the discussion in
\cite[Theorem 2.4]{lop99}. For an approach using variational methods
we refer to \cite[Theorem 4.1]{eg00}).

\section{Convergence to the nonnegative equilibrium, in the case of a bounded domain.}
Theorem \ref{main1} establishes for any $\lambda>\lambda_1$,  the
existence of a unique nonnegative equilibrium point for the semiflow
$\mathcal{S}(t)$. In the light of Theorem
\ref{main2}, in order to prove convergence of solutions of (\ref{eq1.0}) to the nonnegative equilibrium, it remains to verify (a) that solutions of  (\ref{eq1.0}) remain positive for all times and (b) the asymptotic stability of the nonnegative equilibrium.
\begin{proposition}
\label{PositCone}
Assume that condition $(\mathcal{H_{\alpha}})$ or
$(\mathcal{H^{\infty}_{\beta}})$ holds. The set
$$\mathcal{D}_+:=\left\{\phi\in D^{1,2}_0(\Omega,\sigma)\,:\,\phi(x)\geq 0\;\;\mbox{on}\;\; \overline{\Omega}\right\},$$
is a positively invariant set for the semiflow $\mathcal{S}(t)$.
\end{proposition}
{\bf Proof:}\ \ From Proposition \ref{locds}, we have that
solutions are globally defined in time.  It suffices to show that
a kind of  maximum principle holds, that is, solutions of
(\ref{eq1.0}) corresponding to nonnegative initial data, remain
positive.  We adapt an argument from \cite[Proposition
5.3.1]{cazh}. Let  $\phi_0\in D^{1,2}_0(\Omega,\sigma)$, $\phi_0\geq 0$ a.e
in $\Omega$,  and
$\phi\in\mathrm{C}([0,+\infty);{D}^{1,2}_0(\Omega,\sigma))\cap
\mathrm{C}^1([0,+\infty);L^2(\Omega))$ the global in time
solution of (\ref{eq1.0}), with initial condition $\phi_0$. We
consider $\phi^+:=\max\{\phi,0\}$, $\phi^{-}:=-\min\{\phi,0\}$.
Both $\phi^{+}$ and $\phi^{-}$ are nonnegative,
$\phi^{-},\phi^{+}\in
\phi\in\mathrm{C}([0,+\infty);{D}^{1,2}_0(\Omega,\sigma))\cap
\mathrm{C}^1([0,+\infty);L^2(\Omega))$, and we set $\phi=\phi^+-\phi^{-}$.  We get from
(\ref{eq1.0}), that $\phi^{-}$ satisfies the equation
\begin{eqnarray} \label{eqminus}
\partial_t \phi^{-} - \mathrm{div}(\sigma(x)\, \nabla \phi^{-} ) - \lambda\, \phi^{-}+|\phi|^{2\gamma}\phi^{-}=0.
\end{eqnarray}
Multiplying (\ref{eqminus}) by $\phi^{-}$ and integrating over $\Omega$ we obtain
\begin{eqnarray}
\label{enegminus} \frac{1}{2}\frac{d}{dt}||\phi^{-}||^2_{L^2}+
\int_{\Omega}\sigma(x)|\nabla\phi^{-}|^2dx-\lambda||\phi^{-}||^2_{L^2},
+\int_{\Omega}|\phi|^{2\gamma}|\phi^{-}|^2dx=0,
\end{eqnarray}
which implies that
\begin{eqnarray*}
\frac{1}{2}\frac{d}{dt}||\phi^{-}||^2_{L^2}\leq c\,||\phi^{-}||^2_{L^2}.
\end{eqnarray*}
Thus, by Gronwall's Lemma we obtain
\begin{eqnarray}
||\phi^{-}(t)||^2_{L^2}\leq e^{ct}||\phi_0^{-}||^2_{L^2}=0,\;\;\mbox{for every}\;\;t\in [0,+\infty),
\end{eqnarray}
hence $\phi\geq 0$ for all $t\in (0,+\infty)$, a.e. in $\Omega$.\ \ $\diamond$
\setcounter{lemma}{1}
\begin{lemma}
\label{Garding}
Let condition $(\mathcal{H_{\alpha}})$ be fulfilled. The unique nonnegative equilibrium point which exists for $\lambda>\lambda_1$ is uniformly asymptotically stable.
\end{lemma}
We discuss first the stability properties of the zero solution. The linearization about the zero solution which is an equilibrium point for any $\lambda$ is
\begin{eqnarray*}
\label{linear0}
\partial_t\psi-\mathrm{div}(\sigma(x)\nabla\psi)&-&\lambda\psi=0,\;\;x\in\Omega,\nonumber\\
\psi|_{\partial\Omega}&=&0.
\end{eqnarray*}
It follows from (\ref{eigen1}),
that $\phi=0$ is asymptotically stable in $D^{1,2}_0(\Omega,\sigma)$ if $\lambda <\lambda_1$, and unstable in  $D^{1,2}_0(\Omega,\sigma)$ if $\lambda >\lambda_1$.

The linearization around the nonnegative equilibrium point $u$ of (\ref{eq1.0}), is given by
\begin{eqnarray}
\label{linear}
-\mathrm{div}(\sigma (x)\nabla\psi)&-&\lambda\psi+(2\gamma+1)|u|^{2\gamma}\psi=0,\\
\psi |_{\partial\Omega}&=&0,\nonumber
\end{eqnarray}
and we shall see that for the corresponding eigenvalue problem
\begin{eqnarray}
\label{Gard}
-\mathrm{div}(\sigma (x)\nabla\psi)&-&\lambda\psi+(2\gamma+1)|u|^{2\gamma}\psi=\mu\psi,\\
\psi |_{\partial\Omega}&=&0,\nonumber
\end{eqnarray}
zero is not an eigenvalue. The weak formulation of (\ref{Gard}) is
\begin{eqnarray}
\label{weakgard}
A(\psi,\omega)_{\sigma}&:=&\int_{\Omega}\sigma(x)\nabla\psi\nabla\omega\, dx-\lambda\int_{\Omega}\psi\omega\,dx\nonumber\\
&+&(2\gamma+1)\int_{\Omega}|u|^{2\gamma}\psi\omega\,dx=
\mu\int_{\Omega}\psi\omega\,dx,
\end{eqnarray}
for every $\omega\in {D}^{1,2}_0(\Omega,\sigma)$. The symmetric bilinear form $A_{\sigma}:{D}^{1,2}_0(\Omega,\sigma)\times{D}^{1,2}_0(\Omega,\sigma)\rightarrow\mathbb{R}$ defines a Garding form \cite[pg. 366]{zei85}, since
\begin{eqnarray*}
A_{\sigma}(\psi,\psi)\geq ||\psi||^2_{D^{1,2}_0(\Omega,\sigma)}-\lambda||\psi||^2_{\mathrm{L^2(\Omega)}}.
\end{eqnarray*}
Hence, Garding's inequality is satisfied. Then it follows from
Lemmas \ref{lemma2.1}-\ref{lemma2.2} and \cite[Theorem 22.G pg.
369-370]{zei85}, that the problem (\ref{Gard}) has infinitely many
eigenvalues of finite multiplicity, and if we count the eigenvalues
according to their multiplicity, then
\begin{eqnarray}
-\lambda<\mu_1\leq\mu_2\leq\cdots,\;\;\mbox{and}\;\;\mu_j\rightarrow\infty\;\;\mbox{as}\;\;j\rightarrow\infty.
\end{eqnarray}
The smallest eigenvalue can be characterized by the minimization
problem
\begin{eqnarray}
\label{minlin1}
\mu_1=\min A_{\sigma}(\psi,\psi),\;\;\psi\in D^{1,2}_0(\Omega,\sigma),\;\;||\psi||_{L^2}=1.
\end{eqnarray}
The $j$-th eigenvalue, can be characterized by the
minimum-maximum principle
\begin{eqnarray}
\label{minlin2}
\mu_j=\min_{M\in\mathcal{L}_j}\max_{\psi\in M}A_{\sigma}(\psi,\psi).
\end{eqnarray}
where $M=\{\psi\in D^{1,2}_0(\Omega,\sigma)\,:\,||\psi||_{L^2}=1\}$ and $\mathcal{L}_j$ denotes the class of all sets $M\cap L$ with $L$ an arbitrary $j$-dimensional linear subspace of $D^{1,2}_0(\Omega,\sigma)$.

By using similar arguments as for the proof of Lemmas
\ref{lemma2.3}-\ref{lemma2.7}, we may see that for (\ref{Gard}),
the (nontrivial) eigenfunction corresponding to the principal
eigenvalue $\mu_1$ is nonnegative, i.e $\psi_1\geq 0$ a.e. on
$\Omega$. Since $\mu_1,\psi_1$ satisfy (\ref{weakgard}) we get by
setting $\omega=u$ that
\begin{eqnarray*}
\int_{\Omega}\sigma(x)\nabla\psi_1\nabla u\,
dx-\lambda\int_{\Omega}\psi_1 u\,dx +(2\gamma
+1)\int_{\Omega}|u|^{2\gamma}\psi_1 u\,dx=
\mu_1\int_{\Omega}\psi_1 u.
\end{eqnarray*}
On the other hand, by setting $v=\psi_1$ to the weak formula (\ref{forstat}) we get
\begin{eqnarray*}
\int_{\Omega}\sigma(x)\nabla\psi_1\nabla u\,
dx-\lambda\int_{\Omega}\psi_1 u\,dx
+\int_{\Omega}|u|^{2\gamma}\psi_1 u\,dx=0.
\end{eqnarray*}
Subtracting these equations, we obtain that
\begin{eqnarray}
2\gamma\int_{\Omega}|u|^{2\gamma}u\psi_1dx=\mu_1\int_{\Omega}u\psi_1dx,
\end{eqnarray}
which implies that $\mu_1>0$. \ \ $\diamond$

{\bf Proof of Corollary \ref{tsibadyoball}}: The positivity
property of Proposition \ref{PositCone} and Theorem
\ref{main1}, imply that the solution $\phi (\cdot,t)$ converges
towards the set of nonnegative solutions of (\ref{eq1.1}) as
$t\rightarrow\infty$, in $D^{1,2}_0(\Omega,\sigma)$. In fact, it
follows from Lemma \ref{Garding}, that in the case of
$\lambda>\lambda_1$,  for any nonnegative initial condition
$\phi_0$, $\omega(\phi_0)=\{u\}$. On the other hand it is not hard
to check, by following the computations leading to
(\ref{abset1})-(\ref{abset2}), that in the case
$\lambda<\lambda_1$, $\mathrm{dist}(\mathcal{S}(t)\mathcal{B},
\{0\})\rightarrow 0$ as $t\rightarrow\infty$, for every bounded
set $\mathcal{B}\subset D^{1,2}_0(\Omega,\sigma)$. In this case,
the global attractor $\mathcal{A}$ is reduced to $\{0\}$.\
$\diamond$ \vspace{.12cm}
\setcounter{remark}{2}
\begin{remark}
\label{minLyap}
(Minimization of the Lyapunov Function)\
We expect naturally, that the nonnegative equilibrium points, minimize the Lyapunov function
(\ref{Liapunov}) \cite{Achm}. Assume that condition $(\mathcal{H_{\alpha}})$ holds.
It is not hard to check that  $\mathcal{J}$\ is a bounded from below functional on\
$D^{1,2}_0(\Omega,\sigma)$.
Assume further that\ $\lambda<\lambda_1$.\ Then, the variational
characterization of\ $\lambda_1$ (\ref{eq2.1})-
(\ref{eq2.2}), implies that
\[
\mathcal{J}(\phi) \geq \frac{1}{2\gamma+2}
||\phi||^{2\gamma+2}_{L^{2\gamma+2}(\Omega)},
\]
for every\ $\phi \in D^{1,2}_0(\Omega,\sigma)$.\ Hence the trivial
solution is the global minimizer of the functional\
$\mathcal{J}$.\ However, for\ $\lambda>\lambda_1$\ the origin is
no longer the global minimizer of the functional: consider the
function\ $tu_1$,\ where\ $u_1$\ is the normalized nonnegative
eigenfunction associated to\ $\lambda_1$\ and\ $t>0$\ is small
enough. Then from (\ref{Liapunov}), we have that
\begin{eqnarray*}
\mathcal{J}(tu_1) &=& \frac{t^2}{2}
||u_1||^{2}_{D^{1,2}_0(\Omega,\sigma)} - \frac{\lambda_1\, t^2}{2}
||u_1||^{2}_{L^{2}(\Omega)} - \frac{(\lambda-\lambda_1)\, t^2}{2}
||u_1||^{2}_{L^{2}(\Omega)} \\ &&+ \frac{t^{2\gamma+2}}{2\gamma+2}
||u_1||^{2\gamma+2}_{L^{2\gamma+2}(\Omega)} \\
&=& - \frac{(\lambda-\lambda_1)\, t^2}{2}
||u_1||^{2}_{L^{2}(\Omega)} + \frac{t^{2\gamma+2}}{2\gamma+2}
||u_1||^{2\gamma+2}_{L^{2\gamma+2}(\Omega)} < 0.
\end{eqnarray*}
The justification of the Palais-Smale condition follows from Lemma
\ref{lemma2.1} (ii). Then, Ekeland's variational principle implies
the existence of nontrivial minimizers for\ $\mathcal{J}$.\ These
minimizers are the solutions which belong to the branch\
$C_{\lambda_1}$.

Actually\ $C_{\lambda_1}$\ is a pitchfork bifurcation of
supercritical type, where
Principle of Exchange of Stability holds.
\end{remark}
\begin{remark}
\label{degex}
(Degeneracy exponent)\ Condition (\ref{cruc1}) can be written as a restriction on the ``degeneracy'' exponent
\begin{eqnarray}
\label{crucdeg}
0<\alpha\leq\frac{2(1-\gamma(N-2))}{2\gamma +1}:=\alpha^*.
\end{eqnarray}
This is a restriction on the ``rate'' of decrease of the diffusion coefficient $\sigma$
near every point $z\in \sigma^{-1}\{0\}$.  Unique (since the nonlinearity defines a Lipschitz map) and global in time  solutions of (\ref{eq1.0}) exist, converging towards a global attractor, if $\sigma (x)$ decreases more slowly than
$|x-z|^{\delta}$,\ $\delta\in (0,\alpha^*]$, near every point $z\in \sigma^{-1}\{0\}$. As an example we mention the case $N=2$ and $\gamma=1$ (cubic nonlinearity) where $\alpha^*=2/3$. However, as it follows from the discussion in \cite{Ball2}, (\ref{crucdeg}) (or (\ref{cruc1})) does not possibly define a critical exponent, concerning the existence of global attractor. In the case $\alpha>\alpha^*$, the dynamics related to (\ref{eq1.0}), could be investigated through the theory of generalised semiflows \cite{BallNS,Ball2}. As in the case of the damped semilinear wave equation examined in \cite{Ball2}, where uniqueness of solutions is not assumed, one could possibly prove the existence of a global attractor in the case $\alpha>\alpha^*$, under the hypothesis that weak solutions of (\ref{eq1.0}) satisfy the corresponding energy equation.

Possibly, an interesting issue could be, the extension of the bifuraction result and convergence to equilibrium, to the complex evolution equation discussed in \cite{kz??}. Writing the stationary problem as a real system for the real part $\phi_1$ and imaginary part $\phi_2$, one could observe, that at least in the case of the Ginzburg-Landau equation with real coefficients, the stationary problem defines {\em a potential system \cite{Pon}},
which admits a positive principal eigenvalue $\lambda_1$.
The eigenvalue $\lambda_1$ could be a
bifurcation point, from which two global branches bifurcate. These
branches could consist of semitrivial solutions (i.e. solutions of the form
$(\phi_1,0)$ or $(0,\phi_2)$). 
\end{remark}
\begin{remark}
\label{CrucComments}
(Lack of compactess)\ Our approach concerning convergence to the equilibrium, which combines the charachterization of the global attractor of  and global bifurcation theory depends heavily on
$(\mathcal{H_{\alpha}})$ and 
$(\mathcal{H^{\infty}_{\beta}})$, ensuring compactness of the linear and nonlinear operators involved to our study, either in the bounded or unbounded domain case. Thus  it is natural to ask if a  relaxation of the aforementioned conditions which  may  give rise to noncompatness,  could allow for a generalization of the results of Sections 3-5.

A starting point, could be the generalization of the result concerning the existence of the global attractor. One could assume conditions $(\mathcal{H_{\alpha}})$ and 
$(\mathcal{H^{\infty}_{\alpha}})$ for some $\alpha\in [0,+\infty)$. As it is noted in \cite[Remark 2.1]{cm00}, if $\sigma\in L^1_{\mathrm{loc}}(\Omega )$, $\Omega\subseteq\mathbb{R}^N$, $N\geq 2$, satisfies $(\mathcal{H}_{\alpha})$ then it satisfies $(\mathcal{H_{\beta}})$ for any $\beta\geq\alpha$ and if $(\mathcal{H^{\infty}_{\alpha}})$ holds, then $(\mathcal{H^{\infty}_{\delta}})$ is valid for any $\delta\in [0,\alpha]$. 

For example, in the case where $\Omega=\mathbb{R}^N$ one can also consider as an energy space, the space $H^1_0(\mathbb{R}^N,\sigma)$, defined as the closure of $C^{\infty}_0(\mathbb{R}^N)$ with respect to the norm 
\begin{eqnarray*}
||\phi||^2_{H,\sigma}=\int_{\mathbb{R}^N}\sigma (x)|\nabla\phi|^2dx+\int_{\mathbb{R}^N}|\phi|^2dx.
\end{eqnarray*}
The embedding $H^1_0(\mathbb{R}^N,\sigma)\subseteq L^2(\mathbb{R}^N)$ although obviously continuous, is not compact, in the case where $(\mathcal{H}_{\alpha})$ and $(\mathcal{H^{\infty}_{\alpha}})$, hold for some $\alpha\in (0,2]$. Recall also from Remark \ref{remacom}, that  $D^{1,2}_{0}(\mathbb{R}^N,\sigma)$ is not compactly embedded in $L^2(\mathbb{R}^N)$,  if $\sigma$ grows less or equal than quadratically at infinity, even in the case where $(\mathcal{H}_{\alpha})$ is satisfied for some $\alpha\in (0,2)$.

For the existence of global attractors for reaction diffusion equations in unbounded domains, representative references include \cite{bab90,Feiresla}( for semilinear and  degenerate (porous medium) parabolic equations considered on weighted Sobolev spaces) and \cite{BWang99}. The latter provides with an effective remedy for the lack of compactness of the Sobolev imbeddings, with respect to the existence of the global attractor for partial differential equations considered in unbounded domains and in the natural phase space. The idea of \cite{BWang99} is based on the approximation of $\mathbb{R}^N$ by a bounded domain and on the derivation of suitable estimates for the approximation error of the norm of solutions, showing that this approximation error is arbitrary small. These estimates allow for the application of the method developed in \cite{Ball2} which makes use of the energy functional associated to the evolution equation (in \cite{BWang99} a reaction diffusion equation): Asymptotic compactness is shown by passing to the limit of the nonlinear term of the energy functional as the error tends to zero, establishing the existence of a global atractor in $L^2(\mathbb{R}^N)$. 

It would be possibly interesting to attempt to apply this method, to the degenerate equation of the form (\ref{eq1.0}) and investigate if new restrictions could arise between degeneracy, nonlinearity and the parameters involved, through the process of the derivation of an appropriate energy functional, and the estimation of the relevant estimation errors of the generalised Sobolev norms.  

On the other hand, as it is already mentioned in the introduction with respect to the convergence to equilibrium, in the unbounded domain case, one has to deal in general, not only with the lack of compactness but also with the possible appearance of infinite distinct translates of a unique rest point. Thus, it could be also interesting to investigate if the analysis of \cite{Busca}, could provide a framework for the generalization of the convergence result for (\ref{eq1.0}), in the unbounded domain case.

\end{remark}
%
\section{Applications of the Global Bifurcation Result to general elliptic equations}
\setcounter{equation}{0}
We conclude, by mentioning some other examples of degenerate elliptic equations  for which, extensions of the results of Section 4, could be investigated.

{\bf A. Semilinear Equations} We consider the semilinear
problem
\begin{equation} \label{eq1.4}
\begin{array}{ll}
-\nabla(\sigma(x) \nabla u) = \lambda\, f(x)\, u - g(\lambda, x, u), \vspace{0.2cm} \\
\;\;\;\;\;\;\;\;\;\;\;\;\; u |_{\partial \Omega} = 0.
\end{array}
\end{equation}
and the corresponding linear eigenvalue problem
\begin{equation} \label{eq1.5}
\begin{array}{ll}
-\nabla(\sigma(x) \nabla u) = \lambda\, f(x)\, u\, ,\, \vspace{0.2cm} \\
\;\;\;\;\;\;\;\;\;\;\;\;\; u |_{\partial \Omega} = 0,
\end{array}
\end{equation}
where\ $\Omega \subseteq \mathbb{R}^N$,\ $N \geq 2$,\ is an
arbitrary domain (bounded or unbounded).  In this case,
the coefficient functions satisfy:
\vspace{.2cm}

$(\mathcal{F})^{\alpha}$\ $f$\ is a smooth function, at least\
$C^{0,\zeta}_{loc}(\Omega) $,\ for some\ $ \zeta \in (0,1) $,\
such that\ $f \in
L^{\frac{2^{*}_{\alpha}}{2^{*}_{\alpha}-2}}(\Omega)$, and there
exists\ $\Omega^{+}_{f} \subset \Omega$, $|\Omega^{+}_{f}|>0$, such that\ $f(x) > 0$ for
all\ $x \in \Omega^{+}_{f}$.\ \vspace{.2cm}

$(\mathcal{G})^{\alpha}_{\gamma}$\ $g$\ is a Carath\'{e}odory
function, i.e.,\ $g(\cdot,x,\cdot)$\ is a continuous function for
a.e.\ $x \in \Omega$\ and\ $g(\lambda,\cdot,u)$\ is measurable for
all\ $(\lambda,u) \in \mathbb{R}^2$. Moreover, there exist nonnegative
functions\ $c(\lambda) \in C(\mathbb{R})$\ and\ $\rho(x) \in
L^{\infty}(\Omega) \cap L
^{\frac{2^{*}_{\alpha}}{2^{*}_{\alpha}-\gamma}}(\Omega)$,\ such
that\ $|g(\lambda,x,u)| \leq c(\lambda)\, \rho(x)\,
|u|^{\gamma}$,\ for all\ $(\lambda,u) \in \mathbb{R}^2$\ and
almost every\ $x \in \mathbb{R}^N$. \vspace{0.1cm}

Depending on the particular properties of the coefficient
functions, the properties of the global branch could be represented by those of the corresponding approximating problems. For some applications, we
refer to \cite{amb97, arc01, gr00, lop99}.
\vspace{0.2cm} \\
{\bf B. Quasilinear Equations}
We consider quasilinear degenerate elliptic equations of
the form
\begin{equation}\label{eq3b.1}
-\nabla (a(x)\, |\nabla u|^{p-2} \nabla u) = \lambda\, b(x)\,
|u|^{p-2}\ u + f(x)\; |u|^{\gamma-1} u
\end{equation}
and
\begin{equation}\label{eq3b.2}
-\nabla (a(x,u)\, |\nabla u|^{p-2} \nabla u) = \lambda\, b(x)\,
|u|^{p-2}\ u + f(x)\; |u|^{\gamma-1} u,
\end{equation}
where\ $\Omega$\ ia a bounded domain of\ $\mathbb{R}^N$,\ $N\geq
2$\ and\ $1<p<N$.\ A possible treatment could be based on the results
of \cite{dkn97}. Assume that there exists a function\ $\nu(x) \geq
0$ in\ $\Omega$ satisfying \vspace{0.2cm}

$(\mathcal{N}_1)$\, \, $\nu (x) = 0$\ or\ $\nu(x)=\infty$\ in a
finite subset\ $Z \subset \bar{\Omega}$,\ \vspace{0.1cm}

$(\mathcal{N}_2)$\, \, $\nu \in L^{1}_{loc} (\Omega)$,\
$\nu^{-\frac{1}{p-1}} \in L^{1}_{loc} (\Omega)$\ and\ $\nu^{-s}
\in L^{1} (\Omega)$,\ for some\ $s \in \biggl ( \max
\{\frac{N}{p}, \frac{1}{p-1}\}, \infty \biggr )$.\ \vspace{0.3cm}
\\
The coefficient functions\ $a,\; b,\; f$\
satisfy: \vspace{0.2cm}

$(\mathcal{A})$\ $a$\ is a smooth function at least
$C^{0,\zeta}_{loc}$, for some\ $ \zeta \in (0,1) $ a.e. in\
$\Omega$, such that
\[
\frac{\nu(x)}{c} \leq a(x) \leq c\; \nu(x),\;\;\; \mbox{for
some}\;\; c>0.
\]

$(\mathcal{B})$\ $b$\ is a nonnegative and smooth function, at
least\ $C^{0,\zeta}_{loc}(\Omega) $,\ for some\ $ \zeta \in (0,1)
$,\ such that\ $b \in L^{\infty}(\Omega)$,\ \vspace{0.1cm}

$(\mathcal{F})$\ $f$\ is a smooth function, at least\
$C^{0,\zeta}_{loc}(\Omega) $,\ for some\ $ \zeta \in (0,1) $,\
such that\ $f \in L^{\frac{p^{*}_{s}}{p^{*}_{s} - (\gamma+1)}}$,\
where\ $p^{*}_{s} := \frac{Nps}{N(s+1)-ps}$\ and
$p<\gamma+1<p^{*}_{s}$.\ \vspace{0.2cm}

We consider the weighted Sobolev space\ $D_{0}^{1,p}(\Omega,\nu)$ endowed with the norm
\[
||u||_{D_{0}^{1,p}(\Omega,\nu)} := \biggl ( \int_{\Omega} \nu(x)\,
|\nabla u|^p  \biggr )^{1/p} < \infty.
\]
The space\ $D_{0}^{1,p}(\Omega,\nu)$\ is a reflexive Banach space,
enjoying the following embeddings:\vspace{0.1cm}

i)\ \ $D_{0}^{1,p}(\Omega,\nu) \hookrightarrow
L^{p^{*}_{s}}(\Omega)$\ continuously for\ $1<p^{*}_{s}<N$,\
\vspace{0.1cm}

ii)\ $D_{0}^{1,p}(\Omega,\nu) \hookrightarrow L^r (\Omega)$\
compactly for any\ $r \in [1,p^{*}_{s})$.\vspace{0.2cm} \\
For further properties of these spaces, we refer to
\cite{dkn97}, as well as for the proof of the following results.
\begin{lemma}
Assume that conditions\ $(\mathcal{N}_2)$,\ $(\mathcal{A})$,\
$(\mathcal{B})$\ and\ $(\mathcal{F})$\ hold. Then, the
corresponding to (\ref{eq3b.1}) eigenvalue problem
\begin{equation} \label{eq3b.3}
-\nabla (a(x)\, |\nabla u|^{p-2} \nabla u) = \lambda\, b(x)\,
|u|^{p-2}\ u,
\end{equation}
admits a positive principal eigenvalue\ $\lambda_1$,\ given by
\[
\lambda_1 = \inf_{\int_{\Omega} b(x)\, |\phi|^p\, dx =1}
\int_{\Omega} a(x)\, |\nabla \phi|^p\; dx.
\]
Moreover,\ $\lambda_1$\ is simple with a nonnegative associated
eigenfunction\ $u_1$. In addition,\
$\lambda_1$\ is the only eigenvalue with
nonnegative associated eigenfunction.
\end{lemma}
\begin{lemma}
Assume that  conditions\ $(\mathcal{N}_2)$,\ $(\mathcal{A})$,\
$(\mathcal{B})$\ and\ $(\mathcal{F})$\ hold. Then, any weak
solution\ $u \in D_{0}^{1,p}(\Omega,\nu)$\ of (\ref{eq3b.1})
belongs to\ $L^{\infty}(\Omega)$.\ Moreover,\ $u \in
C^{0,\zeta}_{loc}$\ for some\ $ \zeta \in (0,1) $,\ a.e. in\
$\Omega$.\
\end{lemma}
\setcounter{proposition}{2}
\begin{proposition}
\label{genDrab}
Assume that  conditions\ $(\mathcal{N}_2)$,\ $(\mathcal{A})$,\
$(\mathcal{B})$\ and\ $(\mathcal{F})$\ hold. Then the principal
eigenvalue\ $\lambda_1$\ of (\ref{eq3b.3}), is a bifurcation point
of the problem (\ref{eq3b.1}).
\end{proposition}
Proposition \ref{genDrab}, could be extended to a global bifurcation result as follows: Assuming that  $\nu$\ satisfies in adition, condition\
$(\mathcal{N}_1)$, the principal eigenvalue\ $\lambda_1$\ is a
bifurcating point of a global branch. We may
adapt the same procedure described in Sections 2-4, by considering
similar approximating problems. It is interesting to note that in this case,
Picone's identity is still applicable.
\setcounter{theorem}{3}
\begin{theorem}
Assume that the conditions\ $(\mathcal{N}_1)$,\
$(\mathcal{N}_2)$,\ $(\mathcal{A})$,\ $(\mathcal{B})$\ and\
$(\mathcal{F})$\ hold. Then the branch bifurcating from the
principal eigenvalue\ $\lambda_1$\ of (\ref{eq3b.3}), is a global
branch of solutions of the problem (\ref{eq3b.1}). Moreover, any
solution which belongs in this branch, is nonnegative.
\end{theorem}

Concerning problem (\ref{eq3b.2}), we assume that\ $a(x,u)$,\
$b(x,u)$, are sufficiently smooth functions satisfying the following
conditions: \vspace{.2cm}
\newline
$(\mathcal{AS})\ $\ $\frac{\nu(x)}{c} \leq a(x,s) \leq c\; g(|s|)\, \nu(x)$,
$0 \leq b(x,s) \leq b(x)$ and $\lim_{s \to 0} a(x,s) = a(x)$, $\lim_{s \to 0} b(x,s) =
b(x)$ uniformly for a.e. $x \in \Omega$.
\vspace{.2cm}
\newline
Here\ $c>0$,\ $g$\ is a nondecreasing bounded function and\ $a,\;
b$\ satisfy conditions\ $(\mathcal{A})$\ and\
$(\mathcal{B})$,\ respectively. \vspace{0.2cm}

Based again on \cite{dkn97},
and the analysis of Sections 2-4, we may prove
\begin{theorem}
Assume that condition $(\mathcal{AS})$ holds. Then the principal
eigenvalue\ $\lambda_1$\ of (\ref{eq3b.3}), is a bifurcation point
of the problem (\ref{eq3b.2}). Moreover, the corresponding branch
is global, and any solution which belongs to this branch, is
nonnegative.
\end{theorem}
The quasilinear problems,
\[
-\nabla ((k\, |x|^{\alpha} + m\, |x|^{\beta} )\, |\nabla u|^{p-2}
\nabla u) = \lambda\, |u|^{p-2}\ u + f(x)\; |u|^{\gamma-1} u
\]
and
\[
-\nabla ((k\, |x|^{\alpha} + m\, |x|^{\beta} )\, (1+e^{-1/u^2})
|\nabla u|^{p-2} \nabla u) = \lambda\, e^{-u^2}\, |u|^{p-2}\ u +
f(x)\; |u|^{\gamma-1} u,
\]
for some\ $0< \alpha < \min\{p,\; N(p-1)\}$,\ $-N < \beta <0$\
and\ $k$,\ $m$\ nonnegative constants, could serve as examples.

%
{\bf Acknowledgments.}\ \ We would like to thank the referee for his/her valuable comments and suggestions, which improved considerably the presentation of the manuscript. Part of this work was done while
N. B. Zographopoulos had a visiting position at
the Department of Applied Mathematics, University of Crete. The research of N. I. Karachalios was
partially sponsored by a grant from IKY-State Scholarship's Foundation of Greece, contract No. 349, for Postdoctoral Research in the
Department of Mathematics-University of the Aegean. The authors acknowledge financial support from project ``PYTHAGORAS-National Technical University of Athens'' under proposal ``Dynamics of infinite dimensional discrete and continuous dynamical systems and applications''. 
%
%
\bibliographystyle{amsplain}

\end{document}